\documentclass[11pt,reqno]{article}
\usepackage[all]{xy}
\usepackage{amsmath}
\usepackage{amsthm}
\usepackage{amscd}
\usepackage{amssymb}

\input{epsf.tex}

\numberwithin{equation}{section}

\newcommand{\centeredepsfbox}[1]{\centerline{\epsfbox{#1}}}
\newcommand{\nb}[1]{#1\nobreakdash-}

\newcommand{\textmatrix}[4]{\bigl( \begin{smallmatrix} #1 & #2 \\ #3 & #4 
\end{smallmatrix} \bigr)}

\theoremstyle{definition}
\newtheorem*{Definition}{Definition}
\newtheorem*{Remark}{Remark}

\theoremstyle{plain}
\newtheorem{theorem}{Theorem}[section]
\newtheorem{proposition}[theorem]{Proposition}
\newtheorem{lemma}[theorem]{Lemma}
\newtheorem{corollary}[theorem]{Corollary}
\newtheorem{claim}[theorem]{Claim}
\newtheorem{conjecture}[theorem]{Conjecture}
\newtheorem{Problem}[theorem]{Problem}

\newcounter{remarks}
\newenvironment{remarks}%
{\paragraph*{Remarks}\smallskip
     \begin{list}{\arabic{remarks}. }{\usecounter{remarks}%
          \setlength{\leftmargin}{0in}%
          \setlength{\rightmargin}{0in}%
          \setlength{\labelsep}{0pt}%
          \setlength{\labelwidth}{0pt}%
          \setlength{\listparindent}{0pt}%
     }
}
{
\end{list}
}

\newcommand\inv\inverse

\newcommand\wreath{{\;{\text wr}\;}}

\DeclareMathOperator{\Ends}{Ends}
\DeclareMathOperator{\Aff}{Aff}
\DeclareMathOperator{\rank}{rk}
\DeclareMathOperator{\Length}{Length}
\DeclareMathOperator{\SL}{SL}
\DeclareMathOperator{\GL}{GL}
\DeclareMathOperator\Max{Max}
\DeclareMathOperator\nbhd{Nbhd}
\DeclareMathOperator\diam{diam}
\DeclareMathOperator\BS{BS}
\DeclareMathOperator\Bilip{Bilip}

\newcommand\R{{\mathbf R}}
\newcommand\reals{\R}
\newcommand\Q{{\mathbf Q}}
\newcommand\hyp{{\mathbf H}}
\newcommand\C{{\mathbf C}}
\newcommand\Z{{\mathbf Z}}

\newcommand\solv{{\scshape solv}}

\newcommand\inject{\hookrightarrow}
\newcommand\Sum{\sum}
\newcommand\infinity{\infty}
\newcommand{\bdy}{\partial}
\newcommand{\from}{\colon}
\newcommand\composed{\circ}
\newcommand\suchthat{\bigm|}
\newcommand\inverse{{-1}}
\newcommand\union{\cup}
\newcommand\Union{\bigcup}
\newcommand\abs[1]{\left| #1 \right|}
\newcommand\subgroup{<}

\newcommand\Id{\text{Id}}
\newcommand\A{\mathcal A}
\newcommand\intersect{\cap}

\newcommand\Svarc{\v{S}varc}
\newcommand\restrict{\bigm|}
\newcommand\semidirect{\rtimes}

\DeclareMathOperator\QI{QI}
\newcommand\cross{\times}

\newcommand\ext{{\rm int}}
\newcommand\Haus{{\mathcal H}}
\newcommand\F{{\cal F}}
\newcommand\G{{\cal G}}
\newcommand\barM{\overline M}
\newcommand\barN{\overline N}
\newcommand\infdim{inf$\delta$im}
\newcommand\coarsecap{\cap_C}
\renewcommand\O{{\rm O}}
\newcommand\M{{\cal M}}
\newcommand\Mobius{Mobi\"us}
\newcommand\Poincare{Poincar\'e}

\begin{document}

\setcounter{tocdepth}{2}

\title{On the asymptotic geometry of\\ abelian-by-cyclic groups\thanks{To
appear in \emph{Acta Mathematica}}}
\author{Benson Farb\thanks{Supported in part by NSF grant DMS 9704640, 
by IHES, and by the Alfred P. Sloan Foundation.}
\ \ and Lee Mosher\thanks{Supported in part by NSF grant
DMS 9504946 and by IHES.}}
\maketitle


\tableofcontents

\section{Introduction} 

Gromov's Polynomial Growth Theorem \cite{Gromov:PolynomialGrowth} states
that the property of having polynomial growth characterizes virtually
nilpotent groups among all finitely generated groups.

Gromov's theorem inspired the more general problem (see, e.g.
\cite{GhysHarpe:afterGromov}) of understanding to what extent the 
asymptotic geometry of a finitely-generated solvable group determines its
algebraic structure---in short, are solvable groups quasi-isometrically
rigid? In general they aren't: very recently A.\ Dioubina
\cite{Dioubina:solvable} has found a solvable group which is
quasi-isometric to a group which is not virtually solvable; these groups
are finitely generated but not finitely presentable. In the opposite
direction, first steps in identifying quasi-isometrically rigid solvable
groups which are not virtually nilpotent were taken for a special class
of examples, the solvable Baumslag-Solitar groups, in
\cite{FarbMosher:BSOne} and
\cite{FarbMosher:BSTwo}. 

The goal of the present paper is to show that a much broader class of
solvable groups, the class of finitely-presented, nonpolycyclic,
abelian-by-cyclic groups, is characterized among all finitely-generated
groups by its quasi-isometry type. We also give a complete quasi-isometry
classification of the groups in this class; such a classification for
nilpotent groups remains a major open question. Motivated by these
results, we offer a conjectural picture of quasi-isometric
classification and rigidity for polycyclic abelian-by-cyclic groups in
\S\ref{section:polycyclic}.

The proofs of these results lead one naturally from a geometry of groups 
problem to the theory of dynamical systems via the asymptotic
behavior of certain flows and their associated foliations.

\subsection{The abelian-by-cyclic group $\Gamma_M$}
\label{section:abcgroups}

A group $\Gamma$ is {\em abelian-by-cyclic} if there is an exact sequence
$$1\rightarrow A\rightarrow \Gamma\rightarrow Z\rightarrow 1
$$ 
where $A$ is an abelian group and $Z$ is an infinite cyclic group. If
$\Gamma$ is finitely generated, then $A$ is a finitely generated module over
the group ring $\Z[Z]$, although $A$ may not be finitely generated as a
group.

By a result of Bieri and Strebel \cite{BieriStrebel:almostfpsolvable}, the
class of finitely presented, torsion-free, 
abelian-by-cyclic groups may be described
in another way. Consider an $n \cross n$ matrix $M$ with integral
entries and $\det M \ne 0$. Let $\Gamma_M$ be the ascending HNN
extension of $\Z^n$ given by the monomorphism $\phi_M$ with matrix
$M$. Then $\Gamma_M$ has a finite presentation
$$
\langle t,a_1,\ldots ,a_n \suchthat [a_i,a_j]=1, ta_it^{-1}=\phi_M(a_i),
i,j=1,\ldots,n\rangle
$$ 
where $\phi_M(a_i)$ is the word $a_1^{m_1}\cdots a_n^{m_n}$ and the
vector $(m_1,\ldots ,m_n)$ is the $i^{\text{th}}$ column of the matrix
$M$. Such groups $\Gamma_M$ are precisely the class of finitely
presented, torsion-free, 
abelian-by-cyclic groups (see \cite{BieriStrebel:almostfpsolvable} for a
proof involving a precursor of the Bieri-Neumann-Strebel invariant, or
\cite{FarbMosher:BSTwo} for a proof using trees). The group $\Gamma_M$ is
polycyclic if and only if $\abs{\det M}=1$; this is easy to see directly,
and also follows from \cite{BieriStrebel:Valuations}. 

\subsection{Statement of results}

The first main theorem in this paper gives a classification of all
finitely-presented, nonpolycyclic, abelian-by-cyclic groups up to
quasi-isometry. It is easy to see that any such group has a
torsion-free subgroup of finite index, and so is commensurable (hence
quasi-isometric) to some $\Gamma_M$. The classification of these groups
is actually quite delicate---the standard quasi-isometry
invariants (ends, growth, isoperimetric inequalities, etc.) do not
distinguish any of these groups from each other, except that the size of
the matrix $M$ can be detected by large scale cohomological invariants of
$\Gamma_M$.

Given $M\in \GL(n,\R)$, the \emph{absolute Jordan form} of $M$ is the
matrix obtained from the Jordan form for $M$ over $\C$ by replacing each
diagonal entry with its absolute value, and rearranging the Jordan blocks
in some canonical order. 

\begin{theorem}[Classification theorem]
\label{theorem:classification}
Let $M_1$ and $M_2$ be integral matrices with $\abs{\det M_i}>1$ for
$i=1,2$. Then $\Gamma_{M_1}$ is quasi-isometric to 
$\Gamma_{M_2}$ if and only if there are positive integers $r_1,r_2$ such 
that $M_1^{r_1}$ and $M_2^{r_2}$ have the same absolute Jordan form.
\end{theorem}

\begin{Remark} 
Theorem \ref{theorem:classification} generalizes the main result of
\cite{FarbMosher:BSOne}, which is the case when $M_1, M_2$ are positive
$1\times 1$ matrices; in that case the result of \cite{FarbMosher:BSOne}
says even more, namely that $\Gamma_{M_1}$ and $\Gamma_{M_2}$ are
quasi-isometric if and only if they are commensurable. When $n \ge 2$,
however, it's not hard to find $n\cross n$ matrices $M_1, M_2$ such that
$\Gamma_{M_1}, \Gamma_{M_2}$ are quasi-isometric but not commensurable.
Polycyclic examples are given in \cite{BridsonGersten}, and the same ideas
may be used to produce nonpolycyclic examples.
\end{Remark}

The following theorem shows that the algebraic property of being
a finitely-presented, nonpolycyclic, abelian-by-cyclic group is in fact a 
large-scale geometric property.

\begin{theorem}[Quasi-isometric rigidity]
\label{theorem:rigidity}
Let $\Gamma=\Gamma_M$ be a finitely presented abelian-by-cyclic group,
determined by an integer $(n \cross n)$ matrix $M$ with $\abs{\det
M}>1$. Let $G$ be any finitely generated group quasi-isometric to
$\Gamma$. Then there is a finite normal subgroup $K \subset G$ such that
$G/K$ is abstractly commensurable to $\Gamma_N$, for some integer $(n
\cross n)$ matrix $N$ with $\abs{\det N}>1$.
\end{theorem}

\begin{Remark}
Theorem \ref{theorem:rigidity} generalizes the main result of 
\cite{FarbMosher:BSTwo}, which covers the case when $M$ is a positive
$1\times 1$ matrix. The latter result was given a new proof in
\cite{MosherSageevWhyte}, and in  \S\ref{section:qirigidity} we
follow the methods of \cite{MosherSageevWhyte} in proving Theorem
\ref{theorem:rigidity}.
\end{Remark}

\begin{Remark}
The ``finitely presented'' hypothesis in Theorem \ref{theorem:rigidity}
cannot be weakened to ``finitely generated''. Dioubina
shows \cite{Dioubina:solvable} that the wreath product $\Z \wreath \Z$, an
abelian-by-cyclic group of the form $\Z[\Z]$-by-$\Z$, is quasi-isometric
to the wreath product $(\Z\oplus F) \wreath \Z$ whenever $F$ is a finite
group. But $(\Z\oplus F)\wreath\Z$ has no nontrivial finite normal
subgroups, and when $F$ is nonabelian it is not abstractly commensurable
to an abelian-by-cyclic group.
\end{Remark}

One of the key technical results used to prove Theorem
\ref{theorem:classification} 
is the following theorem, which we believe is of independent interest.
It describes a rigidity phenomenon for
\nb{1}parameter subgroups of $\GL(n,\R)$ which generalizes work of 
Benardete \cite{Benardete:divergence} (see also
\cite{Witte:foliations}).

A \nb{1}parameter subgroup $M^t$ of $\GL(n,\R)$ determines a \nb{1}parameter
family of quadratic forms $Q_M(t)=(M^{-t})^T(M^{-t})$ on $\R^n$, where
the superscript~$^T$ denotes transpose. Each $Q_M(t)$ determines a
norm $\| \cdot \|_{M,t}$ and a distance function $d_{M,t}$ on $\R^n$.

\bigskip
\noindent
{\bf Theorem \ref{theorem:parameter:rigidity} (One-parameter subgroup
rigidity)} {\it Let $M^t,N^t$ be 1-parameter subgroups of $\GL(n,\R)$,
such that $M=M^1$ and $N=N^1$ have no eigenvalues on the unit circle. If
there exists a bijection $f\from\R^n\to \R^n$ and constants
$K \ge 1,C\ge 0$ such that for each $t\in \R$ and $p,q \in \R^n$
$$-C+\frac{1}{K} \cdot d_{M,t}(p,q)\leq d_{N,t}(f(p),f(q))\leq K \cdot
d_{M,t}(p,q)+C 
$$
then $M$ and $N$ have the same absolute Jordan form.} 
\bigskip

The proof of Theorem \ref{theorem:parameter:rigidity} is given in 
\S\ref{section:oneparameter}, and shows that in fact $f$ is a
homeomorphism with a reasonably high degree of regularity; see
Proposition~\ref{proposition:regularity}. 

\begin{Remark}
The case of Theorem \ref{theorem:parameter:rigidity} when $f$ is the
identity map follows from a theorem of D. Benardete
\cite{Benardete:divergence}. See also D. Witte
\cite{Witte:foliations}. Benardete's theorem determines precisely when two
one-parameter subgroups of $\GL(n,\R)$ diverge, and it applies as well
to matrices with eigenvalues on the unit circle.
\end{Remark}

\subsection{Homogeneous spaces}

Using coarse topological and geometrical methods, we reduce the study 
of quasi-isometries of $\Gamma_M$ to that of a certain Lie group $G_M$. 

After squaring $M$ if necessary, we can assume that $\det M >0$ and that $M$ 
lies on a \nb{1}parameter subgroup $M^t$ of $\GL(n,\R)$. The group
$\Gamma_M$ is a cocompact subgroup of the solvable Lie group
$G_M=\R^n\semidirect_M \R$, where $\R$ acts on $\R^n$ by the $1$-parameter
subgroup $M^t$. The group $\Gamma_M$ is discrete in $G_M$ if and only if
$\det M=1$. See section \ref{section:liegroup} for details.

The groups $G_M$, with their left invariant metrics, give a rich and familiar
collection of examples, including: all real hyperbolic spaces, when $M$ is a
constant times the identity; many negatively curved homogeneous spaces, when
$M$ has all eigenvalues $>1$ in absolute value; and 3-dimensional \solv\
geometry, when $M$ is a $2\times 2$ hyperbolic matrix of determinant $1$.
The negatively curved examples associated to a real diagonal matrix with
all eigenvalues $>1$ were studied by Pansu \cite{Pansu:dimension} (and later
Gromov \cite{Gromov:Asymptotic}), who analyzed their quasi-isometric geometry
using the idea of ``conformal dimension''.

We should mention also the result of Heintze \cite{Heintze} that the
class of connected, negatively curved homogeneous spaces consists
precisely of those spaces of the form $N \semidirect \R$ where $N$ is a
nilpotent Lie group, and the action of $\R$ on the Lie algebra has all
eigenvalues strictly outside the unit circle.

\subsection{Outline of proofs}

After preliminary sections, \S\ref{section:linear:algebra} on Linear
Algebra, and \S\ref{section:liegroup} on The Solvable Lie Group $G_M$,
the proof of Theorem \ref{theorem:classification} can be divided into 3
main parts: \S\ref{section:horizontal},\ref{section:oneparameter} on the
Dynamics of $G_M$; \S\ref{section:coarsetop} on Quasi-Isometries of
$\Gamma_M$ via Coarse Topology; and
\S\ref{section:classification} on Finding the Integers, where the pieces
of the proof are put together. The proof of Theorem \ref{theorem:rigidity}
is contained in \S\ref{section:qirigidity} on Quasi-Isometric Rigidity.
Finally we pose some conjectures and problems in
\S\ref{section:questions}.

\subsection*{\S\ref{section:horizontal},\ref{section:oneparameter}:
Dynamics of $G_M$}

In these two sections we classify the Lie groups $G_M$ up to 
{\em horizontal-respecting} quasi-isometry, that is up to 
quasi-isometries $\phi\from G_M\to G_N$ which take each set of the 
form $\R^m\times \{t\}$ to a set of the form $\R^n\times 
\{h(t)\}$ for some function $h$ called the \emph{induced time change}.

\bigskip
\noindent
{\bf Theorem \ref{theorem:levelset:preserving}$'$ (Horizontal respecting
quasi-isometries: special case)} {\it Let $M,N$ lie on \nb{1}parameter
subgroups $M^t, N^t$ of $\GL(n,\R)$, and suppose that $\det M, \det N > 1$.
If there exists a horizontal respecting quasi-isometry $\phi\from G_M\to
G_N$, then there exist real numbers $r,s>0$ so that $M^r$ and $N^s$ have the
same absolute Jordan form.}
\bigskip

\begin{Remark}
In the special case where $M,N$ are diagonalizable with all eigenvalues
$>1$, this can be extracted from work of Pansu \cite{Pansu:dimension}
\emph{without} the assumption that $\phi$ is horizontal respecting. This
special case was later reconsidered by Gromov
(see \cite{Gromov:Asymptotic} Section 7.C), as an application of his
``\infdim'' invariant. Our statement and proof of Theorem
\ref{theorem:horizontal:respecting} is inspired in part by the ideas of
exponential growth rates built into the \infdim\ invariant (see also
comments after Proposition \ref{proposition:height:rigidity}).
\end{Remark}

In \S\ref{section:horizontal} we give a slightly more general version of
this statement, Theorem~\ref{theorem:horizontal:respecting}.

The proof of Theorem \ref{theorem:horizontal:respecting} uses a certain
dynamical system on $G_M$, the ``vertical flow'' which flows upward at unit
speed along flow lines of the form (point)$\cross\R \subset \R^m
\semidirect_M\R$. When $M$ has no eigenvalues on the unit circle this is a
hyperbolic or Anosov flow, and in general it is a partially hyperbolic flow.
We prove Theorem \ref{theorem:horizontal:respecting} in several steps, using
stronger and stronger dynamical properties of flows in $G_M$.

\paragraph*{Step 1 (Foliations Rigidity, Proposition
\ref{proposition:foliations:preserved})} Using the  Shadowing Lemma from
hyperbolic dynamics we show that $\phi$ coarsely respects three dynamically
defined foliations of $G_M$ and $G_N$: the weak stable, weak unstable, and
center foliations. This, together with a result of Bridson-Gersten that
depends in turn on work of Pansu (see Corollary \ref{corollary:unipotent}),
allows reduction to the case where $M,N$ have no eigenvalues on the unit
circle.

\paragraph*{Step 2 (Time rigidity, Proposition
\ref{proposition:height:rigidity})} We show that the induced time change map of
$\phi$ is actually an \emph{affine} map between the time parameters of $G_M$
and $G_N$. After taking a real power of $N$ and composing with a vertical
translation, we can assume that 
$\phi$ preserves the time parameter, that is $h(t)=t$.

\paragraph*{Step 3 (One-parameter subgroup rigidity, Theorem
\ref{theorem:parameter:rigidity})} From Step 2, $\phi$ induces a
quasi-isometry between corresponding level  sets of the time parameter on
$G_M,G_N$, which reduces the proof to Theorem
\ref{theorem:parameter:rigidity}, One-Parameter Subgroup Rigidity. The latter
theorem is proved in \S\ref{section:oneparameter}, by studying
rigidity properties of certain flags of foliations of $\R^n$ associated to
the absolute Jordan form of $M \in \GL(n,\R)$.

\subsection*{\S\ref{section:coarsetop}: Quasi-Isometries of $\Gamma_M$ via
Coarse Topology} Given an integer matrix $M \in \GL(n,\R)$ with
$\det M >1$, we study the geometry of $\Gamma_M$ by
constructing a  contractible metric cell complex $X_M$ on which
$\Gamma_M$ acts freely, properly discontinuously and cocompactly by
isometries, so that $\Gamma_M$ is quasi-isometric to $X_M$. 
Topologically, $X_M$ is a product of $\R^{m}$ with the 
homogeneous directed tree $T_M$ with one edge entering and $d$ 
edges leaving each vertex. Here $d=\det M$. Metrically, 
for every coherently oriented line $\ell$ in $T_M$, the metric 
on $X_M$ is such that $\R^m\times \ell$ is isometric to $G_M$. 

The main result of this section, Proposition \ref{proposition:induced},
says that a quasi-isometry $f\from X_M\to X_N$ induces a quasi-isometry
$\phi\from G_M\to G_N$ which respects horizontal foliations. This is
proved using coarse geometric and topological methods. This is precisely
where the condition $\det M, \det N > 1$ is essential for our
proof, since it gives that the trees $T_M, T_N$ have nontrivial
branching, and this branching allows us to show that~$f$ ``remembers''
the branch points (see Step 2 of \S\ref{section:coherent}). 

While this proof is in the spirit of \cite{FarbMosher:BSOne}, further
complications arise in this more general case
(see \S\ref{section:coherent}). Also, for other applications (e.g.\
\cite{FarbMosher:sbf},
\cite{MosherSageevWhyte}), we shall derive Proposition
\ref{proposition:induced} from a still more general result, Theorem
\ref{theorem:horizontal}, which applies to many graphs of groups whose
vertex and edge groups are fundamental groups of aspherical manifolds of
fixed dimension.

\subsection*{\S\ref{section:classification}: Finding the integers}
Given integer matrices $M,N \in \GL(n,\R)$ with $\abs{\det M},
\abs{\det N} > 1$ such that $\Gamma_M$ and $\Gamma_N$ are quasi-isometric,
a simple argument allows us to reduce to the case of positive determinant,
and then the results of \S\ref{section:horizontal}--\ref{section:coarsetop} 
combine to show that there are positive real numbers $r,s$ so that $M^r$
and $N^s$ have the same absolute Jordan form. We need to show that {\em
integral} $r,s$ exist. This is done by showing that a quasi-isometry
$X_M\to X_N$ induces a bilipschitz homeomorphism between certain
self-similar Cantor sets attached to $X_M$ and $X_N$. Applying a theorem
of Cooper on bilipschitz types of these Cantor sets allows us to conclude
that $(\det M)^p = (\det N)^q$ for some integers $p,q \ge 1$,
from which the desired conclusion follows.

\subsection*{\S\ref{section:qirigidity}: Quasi-Isometric Rigidity}  To
prove Theorem \ref{theorem:rigidity}, we use the coarse topology results
from \S\ref{section:coarsetop} to show that a group quasi-isometric to
some $\Gamma_M$ admits a quasi-action on a tree of $n$-dimensional Euclidean
spaces. We then use the results of \cite{MosherSageevWhyte} to convert this
quasi-action into a true action on a tree, whose edge and vertex stabilizers
are finitely generated groups quasi-isometric to $\Z^n$. The proof is
completed by invoking well-known quasi-isometry invariants, combined with a
brief study of injective endomorphisms of virtually abelian groups.

\subsection*{\S\ref{section:questions}: Questions} We make some
conjectures concerning possible extensions of this work to the
polycyclic case. Also, we state some problems on the quasi-isometry group
of $\Gamma_M$.

\paragraph{Acknowledgements} We thank Kevin Whyte and Amie Wilkinson for
all their help. We are also grateful to the IHES, where much of this work
was done.



\section{Preliminaries} 

This brief section reviews some basic material; see for example
\cite{GhysHarpe:afterGromov}.

Given $K\geq 1,C\geq 0$, a $(K,C)$ \emph{quasi-isometry} between metric
spaces is a map $f\from X\to Y$ such that: 
\begin{enumerate}
\item For all $x, x' \in X$ we have
$$\frac{1}{K} \cdot d_X(x,x') - C \le d_Y(f(x),f(x')) \le K \cdot
d_X(x,x') + C
$$
\item For all $y\in Y$ we have $d_Y(y,f(X))\leq C$.
\end{enumerate}
If $f$ satisfies (1) but not necessarily (2) then it is called a 
$(K,C)$ \emph{quasi-isometric embedding}. If $f$ satisfies only the
right hand inequality of (1) then $f$ is $(K,C)$ \emph{coarsely
lipschitz}, and if in addition $C=0$ then $f$ is
\nb{$K$}\emph{lipschitz}. 

A \emph{coarse inverse} of a quasi-isometry $f \from X \to Y$ is a
quasi-isometry $g \from Y \to X$ such that, for some constant $C'>0$, we
have $d(g\circ f(x),x)<C'$ and $d(f\circ g(y),y)<C'$ for all $x \in X$
and $y \in Y$. Every $(K,C)$ quasi-isometry $f \from X\to Y$ has a
$K,C'$ coarse inverse $g \from Y \to X$, where $C'$ depends only on
$K,C$: for each $y\in Y$ define $g(y)$ to be any point $x \in X$ such
that $d(f(x),y) \le C$.

A fundamental fact observed by Efremovich, by Milnor
\cite{Milnor:curvature}, and by \Svarc, which we use repeatedly without
mentioning, states that if a group $G$ acts properly discontinuously and
cocompactly by isometries on a proper geodesic metric space $X$, then $G$
is finitely generated, and $X$ is quasi-isometric to $G$ equipped with
the word metric. 

Given a metric space $X$ and $A,B \subset X$, we denote the
\emph{Hausdorff distance} by 
$$d_\Haus(A,B) = \inf \{r \in [0,\infinity] \suchthat A \subset N_r(B)
\quad\text{and}\quad B \subset N_r(A)\}
$$

The following lemma says that an ambient quasi-isometry induces a
quasi-isometry between subspaces of a certain type. A map $\sigma \from
S \to X$ between geodesic metric spaces is \emph{uniformly proper} if
there is function $\rho \from [0,\infinity) \to [0,\infinity)$ with
$\displaystyle\lim_{t \to \infinity} \rho(t)=+\infinity$, and
constants $K \ge 1, C \ge 0$, such that for all $x,y \in S$ we have:
$$\rho\bigl(d_S(x,y)\bigr) \le d_X(\sigma(x),\sigma(y)) \le K \cdot
d_S(x,y) + C
$$
The function $\rho$ and the constants $K,C$ are called \emph{uniformity
data} for $\sigma$.

\begin{lemma}
\label{lemma:path:subspace}
Given geodesic metric spaces $X,Y,S,T$, a quasi-isometry $f \from X \to
Y$, and uniformly proper maps $\sigma \from S \to X$ and $\tau\from T\to
Y$, suppose that $d_\Haus\bigl(f\sigma(S),\tau(T)\bigr) < \infinity$.
Then $S,T$ are quasi-isometric. To be explicit, any function
$g\from S \to T$ such that $d_Y\bigl(f\sigma(x),\tau g(x)\bigr)$ is
uniformly bounded is a quasi-isometry; the quasi-isometry constants for
$g$ depend only on those for $f$, the uniformity data for $\sigma$ and
$\tau$, and the bound for $d_Y\bigl(f\sigma(x),g\tau(x)\bigr)$.
\end{lemma}

\begin{proof} Pick $K\ge 1$, $C \ge 0$ and $\rho \from [0,\infinity) \to
[0,\infinity)$ such that $f$ is a $(K,C)$ quasi-isometry,
$d_Y\bigl(f\sigma(x),g\tau(x)\bigr) \le C$, and $\rho,K,C$ are
uniformity data for $\sigma,\tau$. 

Consider $x,y \in S$ such that $d_S(x,y) \le 1$. We have
$d_Y(f\sigma(x),f\sigma(y)) \le K^2 + KC + C$, and so
$d_Y(\tau g(x),\tau g(y)) \le K^2 + KC + 3C$, from which it follows
that $\rho\bigl(d_T(g(x),g(y))\bigr) \le K^2+KC+3C$. Since
$\displaystyle\lim_{t \to \infinity} \rho(t)=\infinity$ we obtain a bound
$d_T(g(x),g(y)) \le A$ depending only on $K,C,\rho$. The usual
``rubber band'' argument, using geodesics in $S$ divided into subsegments
of length $1$ with a terminal subsegment of length $\le 1$, suffices to
prove that $g$ is $(K',C')$ coarsely lipschitz, with $K',C'$ depending
only on $K,C,\rho$. 

For any $\xi \in T$ there is a point $\overline g(\xi) \in S$ such that
$d_Y(f\sigma\overline g(\xi),\tau(\xi)) \le C$. For any $\xi,\eta \in
T$ with $d(\xi,\eta) \le 1$ we have
\begin{align*}
d_Y(f\sigma\bar g(\xi),f \sigma \bar g(\eta)) &\le d_Y(f \sigma \bar
g(\xi),\tau(\xi)) + d_Y(\tau(\xi),\tau(\eta)) + d_Y(f \sigma \bar
g(\eta),\tau(\eta)) \\
&\le K+3C
\end{align*}
and so $\rho\bigl(d_S(\bar g(\xi),\bar g(\eta))\bigr) \le  d_X(\sigma\bar
g(\xi),\sigma \bar g(\eta)) \le K^2 + 4KC$. As above we obtain an upper
bound for $d_S(\bar g(\xi),\bar g(\eta))$ and the rubber band argument
shows that $\overline g$ is coarsely lipschitz. 

For any $x \in S$, setting $\xi=g(x) \in T$, we have
\begin{align*}
d_Y(f\sigma(x),f\sigma\overline g(\xi)) &\le
  d_Y(f\sigma(x),\tau g(x)) + d_Y(\tau(\xi),f\sigma\overline g(\xi)) \\
&\le 2C
\end{align*}
It follows that $d_X(\sigma(x),\sigma\overline g(\xi)) \le 3KC$ and
so $\rho\bigl(d_S(x,\overline g g(x))\bigr) = \rho\bigl(d_S(x,\overline
g(\xi))\bigr)\le 3KC$, yielding an upper bound for $d_S(x,\overline g
g(x))$. Similarly, $d_Y(\xi,g\overline g(\xi))$ is
bounded for all
$\xi \in T$.

Knowing that $g \from S \to T$ and $\overline g \from T \to S$ are coarse
lipschitz maps which are coarse inverses of each other, it easily
follows that $g$ is a quasi-isometry, with quasi-isometry constants
depending only on the coarse lipschitz constants for $g$ and $\overline
g$, and on the coarse inverse constants for $g,\overline g$. 
\end{proof}

\vfill\break


\section{Linear Algebra} 
\label{section:linear:algebra}

In this section we collect some basic results about canonical forms of
matrices, and growth of vectors under the action of a matrix.
 
Let $\M(n,F)$ denote all $n \cross n$ matrices over a field $F$, and let
$\GL(n,F)$ be the group of invertible matrices. Let $\GL_0(n,\R)$ be 
the
identity component of $\GL(n,\R)$, consisting of all matrices of positive
determinant. 

\subsection{Jordan Forms}

A matrix $J \in \M(k,\C)$ is a \emph{Jordan block} it it has the form
$J = J(k,\lambda) = \lambda \cdot \Id + N$ where $\lambda \in \C$ and
$N_{ij} =\delta(i+1,j)$, so $N$ is the $k \cross k$ matrix with $1$'s on
the superdiagonal and $0$'s elsewhere.

A matrix $M \in \M(n,\C)$ is in \emph{Jordan form} if it is
in block diagonal form
$$M = \left( \begin{array}{cccc}
J_1 & 0   & \ldots   & 0   \\
0   & J_2 & \ldots   & 0   \\
\vdots & \vdots & \ddots & \vdots \\
0 & 0 & \vdots & J_I
\end{array}\right)
$$
where each $J_i$ is a Jordan block. Every matrix in $\M(n,\C)$ is
conjugate, via an invertible complex matrix, to a matrix in Jordan
form, unique up to permutation of the Jordan blocks. When all
eigenvalues are real, say $J_i$ has eigenvalue $\ell_i$, we resolve
the nonuniqueness by requiring $\ell_1\ge \ell_2
\ge \cdots \ge \ell_I$, and for each $i=1,\ldots,I-1$, if $\ell_i =
\ell_{i+1}$ then $\rank(J_i) \ge\rank(J_{i+1})$. 

A matrix $J \in \M(k,\R)$ is a \emph{real Jordan block} if it has one of
the following two forms. The first form is an ordinary Jordan block
$J(k,\ell)$ where $\ell \in \R$. The second form, which requires
$k$ to be even, has a $2 \cross 2$ block decomposition of the form
$$J = J(k,a,b) =
\left(\begin{array}{ccccc}
Q(a,b) & \Id & \ldots & 0 & 0 \\
0 & Q(a,b) & \ldots & 0 & 0 \\
\vdots & \vdots & \ddots & \vdots & \vdots \\
0 & 0 & \ldots & Q(a,b) & \Id \\
0 & 0 & \ldots & 0 & Q(a,b)
\end{array}\right)
$$
where $\Id$ is the identity, $0$ is the zero matrix,
$Q(a,b) = \textmatrix{a}{-b}{b}{a}$, and $b \ne 0$. 

A matrix $M \in \M(n,\R)$ is in \emph{real Jordan form} if it is
in block diagonal form as above where each block $J_i$ is a real
Jordan block. Every matrix in $\M(n,\R)$ is conjugate, via an
invertible real matrix, to a matrix in real Jordan form, unique up to
permutation of blocks.

The \emph{absolute Jordan form} of $M \in \M(n,\R)$ is the matrix
obtained from the Jordan form of $M$ by replacing each diagonal entry
$\lambda$ by $\ell=\abs{\lambda}$, and permuting the blocks to resolve the
nonuniqueness. If $M$ is invertible then the absolute Jordan form of $M$
can be written in block diagonal form
$$
\left(\begin{array}{ccc}
J^+_M & 0 & 0 \\
0   & J^0_M & 0 \\
0   &  0  & J^-_M
\end{array}\right)
$$
where the diagonal entries of $J^+_M$ are $>1$, of $J^0_M$ are $=1$, and
of $J^-_M$ are $<1$. We call $J^+_M$ the \emph{expanding part} of the
absolute Jordan form, $J^0_M$ the \emph{unipotent part}, and $J^-_M$ the
\emph{contracting part}, and the block matrix
$\textmatrix{J^+_M}{0}{0}{J^-_M}$ is called the \emph{nonunipotent part}.
Of course, one or more of these parts might be empty. 

Note that the Jordan form of the real matrix $J(k,a,b)$ is
$$
\left(\begin{array}{cc}
J(k/2,a+bi) & 0 \\
0 & J(k/2,a-bi)
\end{array}\right)
$$
and so the absolute Jordan form of $J(k,a,b)$ is
$$
\left(\begin{array}{cc}
J(k/2,\sqrt{a^2+b^2}) & 0 \\
0 & J(k/2,\sqrt{a^2+b^2})
\end{array}\right)
$$
Given $M \in \M(n,\R)$, this process may be applied block-by-block to the 
real Jordan form of $M$, and the blocks then permuted, to obtain the
absolute  Jordan form of $M$.

Let $\GL_\times(n,\R)$ denote the set of all matrices in $\GL(n,\R)$
lying on a \nb{1}parameter subgroup of $\GL(n,\R)$, so $\GL_\times(n,\R)
\subset \GL_0(n,\R)$. It is well-known and easy to see, given a matrix
$M\in\GL(n,\R)$, that $M \in \GL_\times(n,\R)$ if and only if the
negative eigenvalue Jordan blocks of $M$ may be paired up so that the two
blocks occuring in each pair are identical to each other, and this
occurs if and only if $M$ has a square root in $\GL(n,\R)$. Thus, if $M$
does not already lie on a \nb{1}parameter subgroup then $M^2$ does. We
are therefore free to replace a matrix by its square in order to land on a
\nb{1}parameter subgroup.

Given a \nb{1}parameter subgroup $\rho(t)$ of $\GL(n,\R)$, if
$M=\rho(1)$ then we will often abuse notation and write
$\rho(t) = M^t$, despite the fact that $M$ may not lie on a unique
\nb{1}parameter subgroup. 

Given $A \in \M(n,\R)$ in Jordan form---no $J(k,a,b)$ blocks---we say that
$\rho(t) = e^{At}$ is a \nb{1}parameter \emph{Jordan subgroup}. Notice 
that the matrices $e^{At}$ are \emph{not} themselves in Jordan form. For
example when $A = J(n,\ell) = \ell \cdot \Id + N$ is a single $n \cross
n$ Jordan block then $e^{At}$ is obtained by multiplying the scalar
$e^{\ell t}$ with the matrix
\begin{equation}
e^{N \cdot t} = \Sum_{i=0}^n \frac{1}{i!} N^i \cdot t^i =  
\left(
\begin{array}{cccccc}
1 & t & \frac{t^2}{2!} & \cdots & \frac{t^{n-1}}{(n-1)!} & \frac{t^n}{n!} 
\\
  & 1 & t & \cdots & \frac{t^{n-2}}{(n-2)!} & \frac{t^{n-1}}{(n-1)!} \\
  &   & 1 & \cdots & \frac{t^{n-3}}{(n-3)!} & \frac{t^{n-2}}{(n-2)!} \\
  &   &   & \ddots & \vdots & \vdots \\
  &   &   &   & 1 & t \\
  &   &   &   &   & 1
\end{array}
\right) 
\label{equation:matrix}
\end{equation}
Nevertheless, for any Jordan form matrix $J = \ell \cdot \Id + N$ with
$\ell \in \R$, the Jordan form of $e^{J}$ is $e^\ell \cdot \Id + N$.

Given a general \nb{1}parameter subgroup $e^{\mu t}$ in $\GL(n,\R)$, 
choose $A$
so that $A^\inverse \mu A$ is in real Jordan form, and so $A^\inverse \mu 
A =
\delta + \nu + \eta$ where $\delta$ is diagonal, $\nu$ is superdiagonal, 
and
$\eta$ is skew-symmetric. We then have
$$e^{\mu t} = (A e^{(\delta + \nu)t} A^\inverse) (A e^{\eta t}
A^\inverse)
$$
Since $\eta$ is skew symmetric it follows that $e^{\eta t}$ is in the
orthogonal group $\O(n,\R)$. We have therefore proved (see
\cite{Witte:foliations} for this particular formulation):

\begin{proposition}[1-parameter Real Jordan Form]
\label{prop:linear:algebra}
Let $M^t$ be a \nb{1}parameter subgroup of $\GL(n,\R)$. There exists a
\nb{1}parameter Jordan subgroup $e^{Jt}$, a matrix $A\in\GL(n,\R)$, and a
bounded \nb{1}parameter subgroup $P^t$ conjugate into the orthogonal group
$\O(n,\R)$, such that $e^J$ is the absolute Jordan form of
$M$, and letting $\barM^t = A^\inverse e^{Jt} A$ we have
$$M^t = \barM^t P^t = P^t \barM^t
$$
\end{proposition}

\begin{Remark}
In \cite{Witte:foliations} the subgroup $\barM^t$ is called the
\emph{nonelliptic part} of $M^t$, and $P^t$ is called the
\emph{elliptic part}. These two \nb{1}parameter subgroups, which
commute with each other, are uniquely determined by $M^t$.
\end{Remark}

\subsection{Growth of vectors under a linear transformation}

Consider a \nb{1}parameter subgroup $M^t$ of
$\GL(n,\R)$, with real Jordan form $M^t = (A^\inverse e^{Jt} A) P^t =
\barM^t P^t$. Let
$$0 < \lambda_1 < \cdots < \lambda_L
$$
be the eigenvalues of $\barM$. Let $V_l = \ker((\lambda_l \cdot \Id -
\barM)^m)$ be the \emph{root space} of the eigenvalue $\lambda_l$,
where $m$ is the multiplicity of $\lambda_l$. Let $n_l$ be the
\emph{index of nilpotency} of $\barM \restrict V_l$, the
smallest integer such that $V_l = 
\ker((\lambda_l \cdot \Id - \barM)^{n_l})$. 
For $i=0,\ldots,n_l-1$ let $V_{l,i} = \ker((\lambda_l \cdot
\Id -\barM)^{i+1})$, so $V_{l,0}$ is the eigenspace of $\lambda_l$ and
$V_{l,n_l-1} = V_l$.  We thus have the \emph{Jordan decomposition} of
$\barM$, which consists of the direct sum of root spaces
$$\R^n = V_1 \oplus \cdots \oplus V_L
$$
together with the \emph{Jordan filtrations}
$$V_{l,0} \subset \cdots \subset V_{l,n_l-1} = V_l, \quad l=1,\ldots,L
$$
This decomposition is uniquely determined by $\barM$, and hence by $M$.

\begin{proposition}[Growth of vectors]
\label{prop:growth}
With the above notation, there exist constants $A,B>0$ with the
following properties. Given $l=1,\ldots,L$ with $\lambda_l \ge 1$, we have:
\begin{description}
\item[Exponential Lower Bound] If $v \in V_l$ and $t \ge 0$ then
$$\| M^t v \| \ge A \lambda_l^t \| v \|
$$
In fact the same lower bound holds if $v \in V_l \oplus V_{l+1} \oplus
\cdots \oplus V_L$.
\item[Exponential$\cdot$Polynomial Upper Bound] 
Given $i=0,\ldots,n_l-1$, if $v \in V_{l,i}$ and $t \ge 1$ then
$$\| M^t v \| \le B \lambda_l^t t^i \| v \|
$$
In fact the same upper bound holds if $v \in (V_1 \oplus \cdots \oplus
V_{l-1}) \oplus V_{l,i}$.
\item[Exponential$\cdot$Polynomial Lower Bound]  
Given $i=0,\ldots,n_l-1$, if $v \in V_{l,i}\setminus V_{l,i-1}$ then there
exists $C_v > 0$ such that if $t \ge 1$ then
$$\| M^t v \| \ge C_v \lambda_l^t t^i
$$
\end{description}
Given $l=1,\ldots,L$ with $\lambda_l \le 0$, similar statements are
true with negative values of $t$.
\end{proposition}

\begin{proof}
We start with the case when $M^t = e^{Jt}$ is a \nb{1}parameter Jordan
subgroup, and the proposition follows by examining each Jordan block \eqref{equation:matrix}.

The second case we consider is when $M^t$ has all positive real
eigenvalues. By Proposition \ref{prop:linear:algebra} we have
$M^t = A^\inverse e^{Jt} A$, and Proposition \ref{prop:growth} follows
immediately from the first case applied to $e^{Jt}$, together with the
fact that $A$ takes the Jordan decomposition of $M^t$ to the Jordan
decomposition of $e^{Jt}$.

In the general case, applying Proposition \ref{prop:linear:algebra} we
have $M^t = (A^\inverse e^{Jt} A) P^t = \barM^t P^t$. We can the apply
the second case to $\barM^t = A^\inverse e^{Jt} A$. Since $P^t$
commutes with $\barM^t$ it follows that $P^t$ preserves the Jordan
decomposition of $\barM^t$. Proposition \ref{prop:growth} then follows
from the boundedness of $P^t$.
\end{proof}


\section{The Solvable Lie Group $G_M$}
\label{section:liegroup}
Recall that $\GL_\cross(n,\R)$ denotes those matrices in
$\GL(n,\R)$ which lie on a \nb{1}parameter subgroup of $\GL(n,\R)$. Also,
each matrix in $\GL_\cross(n,\R)$ has positive determinant. 

Given a matrix $M \in\GL_\cross(n,\R)$ lying on a \nb{1}parameter subgroup
$M^t$ of $\GL(n,\R)$, we associate a solvable Lie group denoted
$G_M$. This is the semidirect product $G_M=\R^n\semidirect_M\R$  with
multiplication defined by
$$(x,t) \cdot (y,s) = (x + M^t y , t+s)
$$
for all $(x,t), (y,s) \in \R^n \cross \R$. We will often identify $G_M = \R^n
\semidirect_M \R$ with the underlying set $\R^n \cross \R$.

\begin{Remark} Although the Lie group $G_M$ depends on more than just the
matrix $M=M^1$ itself---it depends on the entire \nb{1}parameter subgroup
$M^t$---we suppress this dependence in our notation $G_M = \R^n
\semidirect_M \R$. This is justified by the fact that the quasi-isometry
type of $G_M$ depends only on $M$, not on the \nb{1}parameter subgroup
containing $M$ (see the remark after
Proposition~\ref{proposition:choices}). Henceforth, when we say something
like ``Given $M \in \GL_\cross(n,\R)\ldots$'', we will either implicitly
or explicitly choose a \nb{1}parameter subgroup 
$M^t \subgroup \GL(n,\R)$ with $M^1=M$, which in turn determines $G_M$.
\end{Remark}

If $M$ has integer entries then there is a homomorphism $\Gamma_M \to 
G_M$ taking the commuting generators $a_1,\ldots,a_n$ to the standard 
basis of the integer lattice $\Z^n \cross 0 \subset \R^n \cross 0 \subset 
\R^n \cross \R$, and taking the stable letter $t$ to the generator $(0,1) 
\in \R^n \cross \R$. The relator $t a_i t^\inverse = \phi_M(a_i)$ is 
checked by noting that
$$(0,1) \cdot (x,0) \cdot (0,-1) = (Mx,0), \quad\text{for all}\quad x \in 
\R^n
$$
Cocompactness of the image of this homomorphism is obvious. To see that
$\Gamma_M$ embeds in $G_M$ one checks that in the abelian-by-cyclic extension
$1 \to A \to \Gamma_M \to \Z \to 1$, the group $A$ is identified with the
nested union $\Z^n \union M^\inverse(\Z^n) \union M^{-2}(\Z^n) \union \cdots$, in $\R^n$. This also shows that discreteness of $\Gamma_M$ in $G_M$ is
equivalent to $\det M=1$, which is equivalent to $\Z^n = M(\Z^n)$.

For the next several sections we will investigate the geometry of the 
solvable Lie group $G_M$. In this section we begin by showing that $G_M$
and $G_N$  are quasi-isometric if $M,N$ have powers with the same
absolute Jordan form.  Later in \S \ref{section:coarsetop} we will see
that when $M$ has integer  entries, much of the geometry of $\Gamma_M$ is
reflected in the geometry of $G_M$.

We endow $G_M$ with the left invariant metric determined by taking the 
standard Euclidean metric at the identity of $G_M \approx \R^n \cross \R
= \R^{n+1}$. At a point $(x,t) \in \R^n \cross \R \approx G_M$, the
tangent space is identified with $\R^n \cross \R$, and the Riemannian
metric is given by the symmetric matrix
$$
\left(
\begin{array}{ll}
Q_M(t) &0\\
0&1
\end{array}
\right)
$$
where $Q_M(t) = (M^{-t})^T M^{-t}$. For each $t \in \R$, the identification
$\R^n \approx \R^n \cross t \subset G_M$ induces in $\R^n$ the metric
determined by the quadratic form $Q_M(t)$. This metric has distance formula
$$d_{M,t}(x,y) = \| M^{-t}(x-y)\|
$$

\begin{remarks}
\item When $M$ is a $1 \cross 1$ matrix with entry $a>1$, the
group $G_M$ is isomorphic to $\Aff(\R)$, the group of affine
transformations of $\R$, and as a Riemannian manifold $G_M$ is isometric
to a scaled copy of the hyperbolic plane with constant sectional
curvature depending on $a$. 
\item The eigenvalues of $M$ are greater than $1$ in absolute value if 
and only if all sectional curvatures of $G_M$ are negative (see 
\cite{Heintze}). 
\end{remarks}

\begin{proposition}[How the metric on $G_M$ depends on choices]
\label{proposition:choices}
Given \nb{1}parameter subgroups $M^t, N^t$ in $\GL(n,\R)$, suppose
there exist real numbers $r,s>0$ such that $M^r$ and $N^s$ have the
same absolute Jordan form. Then the metric spaces $G_M$ and $G_N$ are
quasi-isometric. To be explicit there exists $A \in \GL(n,\R)$ and $K \ge
1$ such that for each $t \in \R$, the map $v \mapsto A(v)$ is a
$K$-bilipschitz homeomorphism from the metric $d_{M,t}$ to the metric
$d_{N,\frac{s}{r} \cdot t}$; it follows that the map from $G_M = \R^n
\semidirect_M \R$ to $G_N = \R^n \semidirect_N \R$ given by 
$$(x,t) \mapsto \left( Ax,\frac{s}{r} \cdot t\right)
$$
is a bilipschitz homeomorphism from $G_M$ to $G_N$, with bilipschitz
constant $\sup\{K,\frac{s}{r},\frac{r}{s}\}$. 
\end{proposition}

\begin{Remark} The absolute Jordan form of $M^r$ is uniquely determined
by $M$ and $r$: it is the $r^{\text{th}}$ power of the absolute Jordan
form of $M$. It follows in particular that the quasi-isometry type of
$G_M$ depends only on the matrix $M=M^1$, not on the choice of
\nb{1}parameter subgroup $M^t$.
\end{Remark}

\begin{proof}[Proof of Proposition \ref{proposition:choices}]
We proceed in cases.

\medskip
\textbf{Case 1:} 
Assume that $N^t=e^{Jt}$ is the unique \nb{1}parameter Jordan subgroup
such that $N=e^J$ is conjugate to the absolute Jordan form of
$M$. Applying Proposition \ref{prop:linear:algebra} we have
$$M^t = (A^\inverse N^t A) P^t
$$
where $A \in \GL(n,R)$ and the \nb{1}parameter subgroup $P^t$ is bounded. 

Choose $t \in \R$ and $v \in \R^n$. We must show that the two numbers
$\|M^{-t}v\| = \| P^{-t}(A^\inverse N^{-t} A) v\|$ and $\| N^{-t} A
v\|$ have ratio bounded away from $0$ and $\infinity$, with bound
independent of $t,v$. Setting $u = N^{-t}Av$, it suffices to show that
$\|P^{-t} A^\inverse u\|$ and $\|u\|$ have bounded ratio. But this is
clearly true, with a bound of
$$\bigl(\sup_{t} \|P^t\| \bigr) \cdot \Max\{\|A\|,\frac{1}{\|A\|}\}
$$
since the \nb{1}parameter subgoup $P^t$ is bounded.

\medskip
\textbf{Case 2:} 
Assume that there exists $a>0$ such that $M^t = N^{at}$ for all
$t$. Then the metrics $d_{M,t}$ and $d_{N,at}$ are identical.

\medskip
\textbf{General case:} 
Applying Case 2 we may assume that $\det M = \det N$. Applying Case 1
twice we may go from $G_M$ to $G_{e^J}$ to $G_N$, where $e^J$ is
conjugate to the absolute Jordan form of $M$ and of $N$.
\end{proof}


\section{Dynamics of $G_M$, Part I:\\ Horizontal Respecting
Quasi-isometries}
\label{section:horizontal}

In this section we begin studying the asymptotic geometry of the solvable 
Lie groups $G_M$ associated to \nb{1}parameter subgroups $M^t$ of 
$\GL(n,\R)$. As we saw in Section \ref{section:liegroup}, the
quasi-isometry type of $G_M$ depends only on $M$, not on the choice of
\nb{1}-parameter subgroup $M^t$ passing through $M$; see the remark after
Proposition \ref{proposition:choices}. We therefore continue to suppress
the choice of \nb{1}parameter subgroup in our notation.  Further, we do
not restrict the determinant to be $>1$: the results of this section hold
even when $\det M = 1$.

\subsection{Theorem \ref{theorem:horizontal:respecting} on horizontal
respecting quasi-isometries}
\label{section:horizontaldefs}

Let $X,Y$ be metric spaces. Let $\F$ be a decomposition of $X$,
that is, a collection of disjoint subsets of $X$ whose union is
$X$. Let $\G$ be a decomposition of $Y$. Motivated by a foliation
of a manifold, the elements of these decompositions are called
\emph{leaves} and the decomposition itself is called the
\emph{leaf space}.  A quasi-isometry $\phi \from X \to Y$ is said
to \emph{coarsely respect} the decompositions $\F,\G$ if there
exists a number $A \ge 0$ and a map of leaf spaces $h \from \F \to
\G$ such that for each leaf $L \in \F$ we have
$$d_\Haus(\phi(L),h(L)) \le A
$$

For example, consider the space $G_M$. The coordinate function
$G_M \approx \R^n\cross \R \to \R$ given by $(x,t) \mapsto t$ is
called the \emph{time function} of $G_M$. The level sets $P_t
\approx \R^n \cross t$ form the \emph{horizontal foliation} of
$G_M$, whose leaves are called \emph{horizontal leaves} of $G_M$,
and whose leaf space is $\R$. Notice that $d_\Haus(P_s,P_t) =
\abs{s-t}$, and so the time function induces an isometry between the 
horizontal leaf space equipped with the Hausdorff metric and $\R$.

Consider another matrix $N \in \GL_\times(n,\R)$, and denote the
horizontal  leaves of $G_N$ by $P'_t$.

\begin{Definition}[Horizontal respecting]
A quasi-isometry $\phi \from G_M\to G_N$ is said to be
\emph{horizontal respecting} if it coarsely respects the
horizontal foliations of $G_M, G_N$. That is, there exists a
function $h \from \R\to \R$ and $A \ge 0$ such that
$d_\Haus(\phi(P^{\vphantom\prime}_{t}),P'_{h(t)}) \le A$ for all $t\in \R$.

The function $h \from \R \to \R$ is called an \emph{induced time change} for $\phi$, with \emph{Hausdorff constant} $A$.
\end{Definition}

If $h,h'$ are two induced time changes for $\phi$ then $\sup_t
\abs{h(t)-h'(t)} \le A+A' < \infinity$, where $A,A'$ are Hausdorff
constants for $h,h'$ respectively. Also, if $h \from \R \to \R$ is an induced
time change for $\phi$ with Hausdorff constant $A$, if $A' \ge 0$,
and if $h' \from \R\to \R$ is any function satisfying $\sup_{t \in
\R} \abs{h(t)-h'(t)} \le A'$, then $h'$ is also an induced time change for $\phi$, with Hausdorff constant $A+A'$.

\begin{lemma}
\label{lemma:height}
For each $K,C,A$ there exists $C'$ such that if $\phi \from G_M \to
G_N$ is a horizontal respecting $(K,C)$ quasi-isometry, and $h \from
\R \to \R$ is an induced time change for $\phi$ with Hausdorff constant $A$, then $h$ is a $(K,C')$ quasi-isometry of $\R$. 
\end{lemma}

\begin{proof}
We have $\abs{h(t)-h(s)} \le d_\Haus(P_{h(t)},P_{h(s)}) + 2A \le
K\abs{t-s} + C + 2A$. The reverse inequality is similar, and so
$h$ is a quasi-isometric embedding. Since $\phi$ is coarsely onto,
an easy argument shows $h$ is coarsely onto.
\end{proof}

A $(K,C')$ quasi-isometry $h \from \R \to \R$ induces a bijection
of the two-point set $\Ends(\R) = \{-\infinity,+\infinity\}$:
given $\eta_1,\eta_2 \in \Ends(\R)$, we have $h(\eta_1)=\eta_2$ if
and only if $h$ takes every sequence that diverges to $\eta_1$ to
a sequence that diverges to $\eta_2$. The following two properties
of $h$ are equivalent:
\begin{enumerate}
\item $h$ induces the identity on $\Ends(\R)$.
\item $h$ is \emph{coarsely increasing}, that is there exists $L>0$ such 
that if $t>s+L$ then $h(t)>h(s)$.
\end{enumerate}
That (2) implies (1) is obvious. The other direction is true with
any $L > 2C'K$, for if there existed $t \ge s+L$ with $h(t) \le
h(s)$, then since $h$ induces the identity on $\Ends(\R)$ there
would exist $t'>t$ such that $\abs{h(s) - h(t')} \le C'$, but also
$\abs{h(s)-h(t')} \ge \abs{s-t'}/K - C' \ge L/K-C' > C'$, a
contradiction.

If $h \from \R \to \R$ is an induced time change of a horizontal respecting 
quasi-isometry $\phi \from G_M \to G_N$, and if $h$ satisfies the 
equivalent properties (1) and (2), then we say that $\phi$ \emph{coarsely 
respects the transverse orientation} of the horizontal foliations.

\bigskip \noindent {\bf Terminology (time vs.\ height): } In some
contexts the vertical parameter which we have been calling
``time'' will also be called {\em height}, as sometimes this
terminology is more suggestive, for 
example in discussing horizontal foliations.
\bigskip

Here is the main result, whose proof will occupy the remainder of this 
section and the next section.

\begin{theorem}[Horizontal respecting quasi-isometries]
\label{theorem:levelset:preserving}
\label{theorem:horizontal:respecting}
Let $\phi\from G_M\to G_N$ be a quasi-isometry which coarsely
respects the transversely oriented horizontal foliations of $G_M$
and $G_N$. Then there exist real numbers $r,s>0$ so that
$M^r$ and $N^s$ have the same absolute Jordan form.
\end{theorem}

Our proof of Theorem \ref{theorem:horizontal:respecting} 
proceeds in steps, following the outline given in the 
introduction.

\subsection{Step 1a: Hyperbolic dynamics and the Shadowing Lemma}
\label{section:shadowing}

The Lie group $G_M$ has a natural flow which fits into the theory of partially
hyperbolic dynamical systems. From the dynamics we find that the flow has
several invariant foliations, the ``weak stable, weak unstable, and center''
foliations. In \S\S\ref{section:shadowing},\ref{section:foliations:rigidity},
by using the Shadowing Lemma 
(\cite{HirschPughShub}, Lemma 7.A.2, page 133), we prove that a
horizontal respecting quasi-isometry $G_M \to G_N$ also respects the
dynamically defined foliations of $G_M$, $G_N$. 

From this result we obtain the
first piece of our rigidity theorem, by showing that expanding, contracting,
and unipotent parts of the absolute Jordan forms of $M$ and $N$ have the same
ranks respectively, and that the unipotent parts are identical.

\subsubsection{Dynamically defined foliations}
Consider a \nb{1}parameter subgroup $M^t \in \GL(n,\R)$, with real
Jordan form $M^t = \barM^t P^t$. Consider the Jordan decomposition of
$\barM$, and group the root spaces according to whether the
corresponding eigenvalue is $<1$, $=1$, or $>1$ (alternatively, a
logarithm which is $<0$, $=0$, or $>0$), to obtain a decomposition
$\R^n = V^- \oplus V^0 \oplus V^+$. 

\begin{Remark}
It might happen that one or two of the factors $V^-$, $V^0$, $V^+$ is
trivial, that is, \nb{0}dimensional, for instance when all eigenvalues
of $M$ lie outside the unit circle.
\end{Remark}

Now consider the Lie group $G_M = \R^n \semidirect_M \R$ determined by
a \nb{1}parameter subgroup $M^t$. Define the \emph{vertical flow} $\Phi$ on
$G_M$ to be 
$$\Phi_t(x,s) = (x,s+t)$$ 
The tangent bundle $TG_M$ has a 
$\Phi$-invariant splitting 
$$TG_M = E^s \oplus E^c \oplus E^u$$ 
\noindent
defined as follows. The tangent
space at each point $x \in G_M$ is identified with $\R^n \oplus \R$,
and we take 
$$E^s_x = V^- \oplus 0, E^c_x = V^0\oplus \R, E^u_x =
V^+ \oplus 0$$

It is evident from the construction that each of the
distributions $E^s \oplus E^c, E^u\oplus E^c, E^c$ is integrable,
tangent to foliations denoted 
$W^s, W^u, W^c$.  We call these foliations the (weak) \emph{stable,
unstable, and center foliations}, respectively.  
The stable and unstable foliations are transverse, and the intersection of any
stable leaf with any unstable leaf is a center leaf.

Applying the Exponential Lower Bound from Proposition
\ref{prop:growth}, there exist constants $A>0$, $\lambda>1$ such that:
\begin{itemize}
\item If $v \in E^u$, then for $t \ge 0$ we have $\| D \Phi_t v \| \ge
A \lambda^t \| v \|$, and for $t \le 0$ we have $\| D \Phi_t v \| \le
\frac{1}{A} \lambda^t \| v \|$.
\item If $v \in E^s$, then for $t \le 0$ we have $\| D \Phi_t v \| \ge
A \lambda^{-t} \| v \|$, and for $t \ge 0$ we have $\| D \Phi_t v \|
\le \frac{1}{A} \lambda^{-t} \| v \|$.
\end{itemize}
Also, applying the Exponential$\cdot$Polynomial Upper Bound from
Proposition~\ref{prop:growth}, there exists $B>0$ and an integer $n
\ge 1$ such that:
\begin{itemize}
\item If $v \in E^c$, then for $\abs{t} \ge 1$ we have $\| D \Phi_t v
\| \le B \abs{t}^n \| v \|$.
\end{itemize}
When we want to emphasize the dependence of the $V$'s and $E$'s on the
\nb{1}parameter subgroup $M^t$, we will append a subscript, e.g.\ $V^+_M$,
$E^s_M$, etc.

\subsubsection{Shadowing Lemma}

Consider a flow $\Phi$ on a metric space $X$. We write $x \cdot t$ as an
abbreviation for $\Phi_t(x)$. Given $\epsilon,T>0$, an
\emph{$(\epsilon,T)$-pseudo-orbit} of $\Phi$ consists of a sequence of flow
segments $(x_i \cdot [0,t_i])$, where the index $i$ runs over an interval
in $\Z$, such that
$d_X(x_i \cdot t_i,x_{i+1}) < \epsilon$ and $t_i > T$ for all $i$.

\begin{lemma}[Shadowing Lemma]
\label{lemma:shadowing}
Consider a \nb{1}parameter subgroup $M^t$ of $\GL(n,\R)$, and let $\Phi$ be the
vertical flow on $G_M$. For every $\epsilon,T>0$ there exists $\delta,
\epsilon', T' > 0$ such that every $(\epsilon,T)$-pseudo-orbit of $\Phi$ is
{\em $\delta$-shadowed} by an $(\epsilon', T')$-pseudo-orbit of $\Phi$
which is contained in some center leaf. That is, if $(x_i \cdot [0,t_i])$
is an $(\epsilon,T)$-pseudo-orbit, then there is an
$(\epsilon',T')$-pseudo-orbit $(y_i \cdot [0,t_i])$ contained in some
center leaf so that $d(x_i \cdot t,y_i \cdot t) < \delta$ for all~$i$ and
all $t \in [0,t_i]$.
\end{lemma}

\begin{proof} By construction, the foliations $W^s$ and $W^u$ are coordinate
foliations in $\R^{n+1}$; this shows that the flow $\Phi$ has a ``global
product structure'' in the language of hyperbolic dynamical systems.  The lemma
now follows the proof of the Shadowing Lemma in \cite{HirschPughShub}, Lemma
7.A.2, page 133.  A direct proof is also easy to work out, and is left to the
reader.
\end{proof}

\subsection{Step 1b: Foliations rigidity}
\label{section:foliations:rigidity}

The Shadowing Lemma implies further rigidity properties of horizontal
respecting quasi-isometries:

\begin{proposition}[Foliations rigidity]
\label{proposition:foliations:preserved}
Suppose $\phi\from G_M\to G_N$ is a quasi-isometry which coarsely
respects the horizontal foliations and their transverse orientations.
Then $\phi$ also coarsely respects the weak unstable foliations
$W^u_M$, $W^u_N$, the weak stable foliations $W^s_M, W^s_N$, and the
center foliations $W^c_M$, $W^c_N$.
In particular 
\begin{itemize}
\item $\dim(V^+_M) = \dim(V^+_N)$
\item $\dim(V^-_M) = \dim(V^-_N)$
\item $\dim(V^0_M) = \dim(V^0_N)$
\end{itemize}
\end{proposition}

\begin{remarks}
\item
In the case where neither $M$ nor $N$ has any eigenvalue on the unit
circle, the center foliations of both $G_M$ and $G_N$ are simply the
foliations by vertical flow lines, and Proposition 
\ref{proposition:foliations:preserved} says that
$\phi$ respects these foliations. But in the general case, it is not
true that $\phi$ always respects the foliations by vertical flow
lines. For a simple counterexample, consider the $1\times 1$ matrix
$M=N=(1)$, which gives $\Gamma_M =\Gamma_N =
\Z^2$.  There exist horizontal respecting quasi-isometries of $\R^2 =
\R \cross \R$ which do not respect the vertical foliation.

\item
If all eigenvalues of $M$ and $N$ are outside the unit circle, then
both $G_M$ and $G_N$ are negatively curved, and the proposition
follows from a standard fact: a quasigeodesic in a negatively curved
space $X$ is Hausdorff close to a geodesic (this was the approach
taken in \cite{FarbMosher:BSOne} in the case of a $1 \cross 1$ matrix
$M$, where $G_M$ is isometric to a scaled copy of the hyperbolic
plane). This ``fact'' is unavailable when $X=G_M$ is not negatively
curved, forcing us to study horizontal respecting quasi-isometries via
the Shadowing Lemma.
\end{remarks}

Before proving Proposition \ref{proposition:foliations:preserved}, 
we use it to obtain some pieces of our
classification theorem. Since $\rank(J^-_M) = \dim(V^-_M)$ etc., we
immediately have:

\begin{corollary}
If there is a quasi-isometry from $G_M$ to $G_N$ which coarsely
respects the transversely oriented horizontal foliations, then
$\rank(J^-_M) = \rank(J^-_N)$, $\rank(J^0_M) = \rank(J^0_N)$, and
$\rank(J^+_M) = \rank(J^+_N)$.\qed
\end{corollary}

We also have:

\begin{corollary}
\label{corollary:unipotent}
The unipotent blocks of the absolute Jordan forms of $M$ and $N$ are
identical.
\end{corollary}

\begin{proof} 
Let $L$ be some center leaf of $G_M$, of dimension $k$. From Proposition \ref{proposition:foliations:preserved} it follows that $\phi(L)$ is Hausdorff close to some center leaf $L'$ of $G_N$,
also of dimension $k$. By composition with nearest point projection 
(which moves points a uniformly bounded amount) we get an
induced map $L \rightarrow L'$. By Lemma \ref{lemma:path:subspace} this
map is a quasi-isometry. By Proposition \ref{proposition:choices}, 
$L$ and $L'$  are quasi-isometric to the
nilpotent Lie groups $\R^{k-1}\semidirect_{J^0_M} \R$ and
$\R^{k-1}\semidirect_{J^0_N}\R$, respectively. As Bridson and Gersten
have shown \cite{BridsonGersten}, Pansu's invariant \cite{Pansu:CC}
may be used to prove that $J^0_M = J^0_N$.
\end{proof}

\begin{proof}[Proof of Proposition \ref{proposition:foliations:preserved}]
We begin with:
\begin{claim}
\label{claim:quasivertical}
For each vertical flow line $\gamma=\Phi_\R(x)$ in
$G_M$, there exists a center leaf $\tau_\gamma$ in $G_N$ such that
$\phi(\gamma)$ is contained in the $\alpha$-neighborhood of
$\tau_\gamma$, where the constant $\alpha>0$ does not depend on
$\gamma$.
\end{claim}

Before proving the claim, we apply it to prove the proposition as
follows. 

Consider any two vertical flow lines $\gamma_1,\gamma_2$ in $G_M$. By the
claim we have that $\phi(\gamma_1)$ and $\phi(\gamma_2)$ lie,
respectively, in bounded neighborhoods of center leaves $\sigma_1$ and
$\sigma_2$ of $G_N$. Since $h(t) \to \pm\infinity$ as $t \to
\pm\infinity$, for each choice of sign $+$ or $-$ the following two
statements are equivalent, and the second statement implies the third:
\begin{enumerate}
\item The distance between the points $\gamma_1 \intersect P_t$ and
$\gamma_2 \intersect P_t$ in $P_t$ stays bounded as $t \to \pm\infinity$.
\item The distance between the points $\phi(\gamma_1) \intersect P_{h(t)}$
and $\phi(\gamma_2) \intersect P_{h(t)}$ in $P_{h(t)}$ stays bounded as $t
\to \pm\infinity$.
\item The Hausdorff distance between the sets $\sigma_1 \intersect
P_{h(t)}$ and $\sigma_2 \intersect P_{h(t)}$ in $P_{h(t)}$ stays bounded
as $t \to \pm\infinity$.
\end{enumerate} 
Using $-$ signs, the first statement is equivalent to saying that
$\gamma_1,\gamma_2$ are contained in the same unstable leaf of
$G_M$, and the third statement is equivalent to saying that
$\sigma_1,\sigma_2$ are contained in the same unstable leaf of $G_N$. It
follows that $\phi$ takes every unstable leaf of $G_M$ into a bounded
neighborhood of an unstable leaf of $G_N$. Applying the same argument to a
coarse inverse $\overline\phi$ of $\phi$ gives the opposite inclusion.
Since $d(\overline\phi \composed\phi,\Id)<\infinity$ it follows that the
image under $\phi$ of any unstable leaf of $G_M$ lies a bounded Hausdorff
distance from an unstable leaf of $G_N$, that is, $\phi$ coarsely
preserves the unstable foliations. A similar argument using $+$ signs
shows that $\phi$ coarsely preserves stable foliations. By taking
intersections of stable and unstable leaves it follows that $\phi$
coarsely preserves center foliations.

The final statements about dimensions follow from the fact that
dimension is a quasi-isometry invariant, for leaves of the foliations
in question; see \cite{Gersten:dimension} or \cite{BlockWeinberger}.

It remains to prove the claim. Applying Lemma \ref{lemma:height}, we have an
induced time change $h \from \R\to \R$ which is a $(K,C')$-quasi-isometry
with Hausdorff constant $A$, where $C'$ depends only on $K,C,A$. Furthermore by
Lemma \ref{lemma:height} and the comments following it, the map $h$ is coarsely
increasing: there exists $L=L(K,C,A)>0$ such that if $t \ge s+L$ then $h(t) >
h(s)$.

We can furthermore increase $L$, depending only on $K,C',A$, so that:
\begin{multline}\label{equation:increasing}
\text{if } t' \ge t+L,\quad x \in P_{t'},\quad y \in P_t,\quad \phi(x) \in P_{s'} \text{ and } \phi(y) \in P_{s},\\ 
\text{then } s' \ge s + 1.
\end{multline}
In fact taking $L > (C'+2A+1)K$ will do, for then we have
\begin{align*}
h(t') &\ge h(t) + (t'-t)/K - C' \ge h(t') + L/K - C' \\
     &\ge h(t) + 2A+1
\end{align*}
and, since $P_{s'}$ is $A$-Hausdorff close to $P_{h(t')}$ and $P_{s}$
is $A$-Hausdorff close to $P_{h(t)}$, it follows that $s' \ge s+1$.

To prove the claim, we first show that $\phi(\gamma)$ is Hausdorff
close to some pseudo-orbit in $G_N$, and then we apply the Shadowing
Lemma to show that the pseudo-orbit lies in a bounded neighborhood of
some center leaf.

To be more precise, fix a point $x_0 \in \gamma$ and consider the sequence $x_i
= \Phi_{i\cdot L}(x_0)$ for $i \in \Z$.  Let $y_i = \phi(x_i)$, and let $s_i$
be such that $y_i \in P_{s_i}$. From \eqref{equation:increasing} it follows
that $s_{i+1} \ge s_i+1$. Let $t_i = s_{i+1}- s_i \ge 1$.

We claim that there exists $\epsilon>0$, depending ultimately only on
$K,C$, so that $\bigl(y_i \cdot [0,t_i]\bigr)$ is an
$(\epsilon,1)$-pseudo-orbit; in other words, $d(y_i \cdot t_i,y_{i+1})$ is
bounded. To see why, first note that 
\begin{align*}
d(y_i \cdot t_i,y_i) &= t_i = s_{i+1}-s_i \\
   &\le 2A + h(L\cdot (i+1)) - h(L\cdot i)  \\
   &\le 2A + KL + C'
\end{align*}
and then
$$d(y_i,y_{i+1}) \le K\cdot d(x_i,x_{i+1}) + C \le KL+C
$$ 
so we may take $\epsilon=2A+2KL+C+C'$.  

Applying the Shadowing Lemma, there exists $\beta, \epsilon', T'$ such
that $\bigl(y_i \cdot [0,t_i]\bigr)$ is $\beta$-Hausdorff close to an
$(\epsilon',T')$-pseudo-orbit $\bigl(y'_i \cdot [0,t_i]\bigr)$ contained in
some center leaf of $G_N$. On the other hand, since every 
point of $\gamma$ is within distance $L$ of some $x_i$, 
it follows that $\phi(\gamma)$ is uniformly
Hausdorff close to $\bigl(y_i \cdot [0,t_i]\bigr)$, and so it is also
uniformly close to the pseudo-orbit $\bigl(y'_i \cdot [0,t_i]\bigr)$.
\end{proof}

\subsection{Step 2: Time rigidity}
\label{section:height:rigidity}
\label{section:time:rigidity}

The main result of this subsection says that a horizontal respecting 
quasi-isometry has an induced time change function which is be affine.

\begin{proposition}[Time rigidity]
\label{proposition:height:rigidity}
Consider the Lie groups $G_M, G_N$ where $M,N\in \GL_\cross(n,\R)$  each
have  an eigenvalue of absolute value greater than~1. Then there
exists $m \in \R_+$ with the following properties. For all $K \ge 1$, $C,
A \ge 0$ there exists $A' \ge 0$ such that if $\phi\from G_M\to G_N$ is a
$(K,C)$ quasi-isometry which coarsely respects horizontal foliations and
their transverse orientations, with an induced time change of Hausdorff
constant $A$, then there exists $b \in \R$ such that $h(t)=mt+b$ is an
induced time change with Hausdorff constant $A'$. In fact, $m$ can be
computed as follows: Let $\alpha$ (resp. $\beta$) be the least eigenvalue
greater than $1$ of the absolute Jordan form of $M$ (resp.\ $N$); the
numbers $\alpha,\beta$ exist by the assumption on  eigenvalues.  Then
$m=\log \alpha/\log \beta$.
\end{proposition}

\begin{remarks}
\item In the case of self-quasi-isometries of $\Aff(R) = G_{(e^1)}
= \hyp^2$ which coarsely respect the horizontal foliation, this
result is part of Proposition 5.3 of \cite{FarbMosher:BSOne},
where the conclusion is that the induced time change is a
translation of $\R$.
\item One of the delicate points in Gromov's development of the
\infdim\ invariant is the rescaling problem discussed at the
beginning of Section 7.$\mathrm C_1$ of \cite{Gromov:Asymptotic}:
rates of exponential growth change when the parameter is rescaled.
Time Rigidity allows us to avoid the rescaling problem
altogether, by showing that the time parameter is ``natural'' with 
respect to quasi-isometries.
\end{remarks}

\begin{proof} This proof will define a sequence of constants which will 
depend on $K,C,A$ and on the matrices $M$ and $N$. We will indicate
the dependence on $K,C,A$ by writing, for example, $C_1 = C_1(K,C,A)$,
but we will suppress the dependence on $M,N$. Although each constant
in the sequence will depend on previous constants in the sequence, by
induction it will ultimately depend only on $K,C,A,M,N$.

\begin{claim} 
\label{claim:growth}
For each fixed time $t_0$, and for each $t \le t_0$, we have
$$h(t) \ge m(t-t_0) + h(t_0) - C_1
$$
for some $C_1 = C_1(K,C,A) \ge 0$. 
\end{claim}

Accepting this claim for the moment, we prove the
proposition.  The idea is simply that the conclusion of the claim,
applied to both $h$ and its coarse inverse $\bar h$, with $t_0\to
+\infty$, implies the proposition.

Let $s$ be a time parameter for $G_N$. Let $\bar\phi \from G_N \to
G_M$ be a coarse inverse for $\phi$, also a quasi-isometry which
coarsely respects the horizontal foliations and their transverse
orientations, and with an induced time change $\bar h(s)$. The
constants for $\bar\phi$ and $\bar h$ depend only on $K,C,A$. The
claim therefore applies as well to $\bar h$ and we obtain, for each
fixed time $s_0$ and each $s \le s_0$,
$$\bar h(s) \ge \frac{1}{m}(s-s_0) + \bar h(s_0) - C_2
$$
for some $C_2=C_2(K,C,A) \ge 0$.

It is clear $\bar h$ is a coarse inverse for $h$, that is:
$$\abs{\bar h(h(t)) - t} \le C_3, \quad \abs{h(\bar h(s))-s} \le C_3
$$
for some $C_3=C_3(K,C,A) \ge 0$. 

Also, by Lemma \ref{lemma:height} and the comments after it, 
the map $h$ is coarsely
increasing: there exists $L = L(K,C,A) \ge 0$ such that if
$t' > t+L$ then $h(t') > h(t)$.

We reverse the inequality in the claim as follows. Fix $t_0$. Let
$s_0=t_0$.  Consider for the moment some $t \le t_0 - L$. Letting
$s=h(t)$ it follows that $s \le s_0$ and so we have
$$
\bar h(h(t)) \ge \frac{1}{m}(h(t)-h(t_0)) + \bar h(h(t_0)) - C_2
$$
But $t+C_3 \ge \bar h(h(t))$ and $\bar h(h(t_0)) \ge t_0 - C_3$ and so we 
obtain
\begin{align*}
t &\ge \frac{1}{m}(h(t) - h(t_0)) + t_0 - (2C_3+C_2) \\
h(t) &\le m(t-t_0) + h(t_0) + m(2 C_3 + C_2) \\
\end{align*}
This has been derived only for $t \le t_0-L$, but for $t_0 - L \le t 
\le t_0$ we obtain a similar inequality with another constant in place of 
$m(2 C_3+C_2)$. Therefore, for all $t \le t_0$ we obtain
$$m(t-t_0) + h(t_0) - C_4 \le h(t) \le m(t-t_0) + h(t_0) + C_4
$$
for some $C_4 = C_4(K,C,A)$. Note that this is true \emph{for all} 
$t_0$, with $C_4$ independent of $t_0$.

In particular, taking $t_0=0$, for all $t \le 0$ we obtain
$$mt + h(0) - C_4 \le h(t) \le mt + h(0) + C_4
$$
Now take any $t_1 \ge 0$, and since $0 \le t_1$ we obtain
$$m(0-t_1) + h(t_1) - C_4 \le h(0) \le m(0-t_1) + h(t_1) + C_4
$$
and so
$$m t_1 + h(0) - C_4 \le h(t_1) \le mt_1 + h(0) + C_4
$$
Taking $b=h(0)$, this proves that $mt+b$ is an induced time change 
for $\phi$, with Hausdorff constant $A'=C_4+A$.
\medskip

Now we turn to the proof of Claim \ref{claim:growth}.  

Let $M^t = \barM^t Q^t$, $N^t = \barN^t Q'{}^t$ be the real Jordan
forms. Let $U$ (resp.\ $U'$) be the root space with eigenvalue $1$ for
$\barM$ (resp.\ $\barN$). Let $W$ (resp.\ $W'$) be the direct sum of
root spaces with eigenvalue $\ge 1$ for $\barM$ (resp.\ $\barN$). 
Recall that $\alpha$ is the smallest eigenvalue $>1$ for $\barM$,
and $\beta$ is the smallest eigenvalue $>1$ for $\barN$. Let $V$ be
the direct sum of $U$ and the eigenspace with eigenvalue $\alpha$ for
$\barM$. We have $U \subset V \subset W$; let $\F(U), \F(V), \F(W)$ be
the corresponding foliations of $G_M \approx \R^n \cross \R$, whose
leaves are parallel to $U \cross \R, V \cross \R, W \cross R$
respectively. We also have $U' \subset W'$; let $\F(U'), \F(W')$ be
the corresponding foliations of $G_N$.

Here is the idea for proving Claim \ref{claim:growth}. Each leaf of
$\F(V)$ is foliated by leaves of $\F(U)$. Because $V$ is the direct
sum of $U$ with the $\alpha$ eigenspace of $\barM$, it follows that as
$t \to -\infinity$ distinct leaves of $\F(U)$ in $\F(V)$ diverge from
each other \emph{exactly} as $\alpha^{-t}$, measured in the time $t$
horizontal plane of $G_M$. This is a consequence of the Exponential Lower Bound
and the Exponential$\cdot$Polynomial Upper Bound in Proposition
\ref{prop:growth}; notice that it is critical here that $V$ not be the
direct sum of $U$ with the $\alpha$ \emph{root space}, for then
Exponential$\cdot$Polynomial Upper Bound would be at best
$\alpha^{-t}$ times some polynomial, which would mess up the following
calculations. Mapping over via the quasi-isometry $\phi \from G_M \to
G_N$, distinct leaves of $\F(U)$ in a single leaf of $\F(V)$ must
(coarsely) map to distinct leaves of $\F(U')$ in a single leaf of
$\F(W')$, which as $s \to -\infinity$ diverge from each other \emph{at
least as fast as} $\beta^{-s}$, by the Exponential Lower Bound. The
time change map $t \mapsto h(t)=s$ therefore cannot grow slower than
$s=\frac{\log \alpha}{\log \beta} \cdot t$, as $t \to -\infinity$.

To make this precise, pick a leaf $L_V$ of $\F(V)$, contained in some
leaf $L_W$ of $\F(W)$. We use the symbol $\gamma$ to denote a general
leaf of $\F(U)$, which we will typically take to be a subset of
$L_V$. By Proposition \ref{proposition:foliations:preserved}, there
exists a leaf $L_{W'}$ of $\F(W')$ such that
$$
d_\Haus(f(L_W),L_{W'}) \le C_5 = C_5(K,C,A)
$$
and for each leaf $\gamma$ of $\F(U)$ there exists a leaf $\gamma'$ of 
$\F(U')$ such that
$$
d_\Haus(f(\gamma),\gamma') \le C_5
$$
Moreover, if $\gamma \subset L_V$ then $\gamma' \subset L_{W'}$, because 
$L_V \subset L_W$ and so $\gamma'$ stays in a bounded neighborhood of 
$L_{W'}$, but any leaf of $\F(U')$ which is not a subset of $L_{W'}$ 
has points which are arbitrarily far from $L_{W'}$.

Let $P_t$ be the horizontal subset of $G_M$ at height $t \in \R$, and
let $d_t$ denote Hausdorff distance in $P_t$ between closed subsets of
$P_t$. Let $P'_s$ be the horizontal subset of $G_N$ at height $s \in
\R$, and let $d'_s$ denote Hausdorff distance in $P'_s$.

Since the Hausdorff distance in $G_N$ between $\phi(P_t)$ and $P'_{h(t)}$ 
is at
most $A$, the vertical projection from $\phi(P_t)$ to $P'_{h(t)}$
induces a quasi-isometry between $P_t$ and $P'_{h(t)}$; the
multiplicative constant of this quasi-isometry is $K$, and its
additive constant depends only on $K,C,A$. It follows that there
exists a ``coarseness constant'' $C_6=C_6(K,C,A)$ so that for any $t$,
and for any $x,y \in P_t$ with $d_t(x,y) \ge C_6$, if $x',y' \in
P'_{h(t)}$ are the vertical projections of $\phi(x),\phi(y)$ then
\begin{equation}\label{equation:bilip}
\frac{1}{2K} d_t(x,y) \le d'_{h(t)}(x',y') \le 2K d_t(x,y)
\end{equation}

To prove Claim \ref{claim:growth}, fix a time $t_0$ and let $s_0 =
h(t_0)$. Let $\gamma_1,\gamma_2$ be two leaves of $\F(U)$ contained in
$L_V$, and let $\gamma'_i$ be the unique leaf of $\F(U')$ within
bounded Hausdorff distance of $\phi(\gamma_i)$; this bound depends
only on $K,C,A$, as shown in Proposition
\ref{proposition:foliations:preserved}.

In $G_M$, apply the Exponential Lower Bound and the
Exponential$\cdot$Polynomial Upper Bound of Proposition
\ref{prop:growth}, and so for all $t \le t_0$ we have
\begin{align*}
A \cdot \alpha^{-t+t_0} d_{t_0}(\gamma_1 \intersect P_{t_0}, 
                                       \gamma_2 \intersect P_{t_0})
  &\le d_{t}(\gamma_1 \intersect P_t, \gamma_2 \intersect P_t) \\
  &\le B \cdot \alpha^{-t+t_0} d_{t_0}(\gamma_1 \intersect P_{t_0}, 
                                       \gamma_2 \intersect P_{t_0})
\end{align*}
where $A,B$ depend only on $G_M$ (note that $t=t_0$ gives $A \le 1 \le
B$).

We want the distance between $\gamma_1$ and $\gamma_2$ in $P_t$ to be greater
than the coarseness constant $C_6$, for each $t \le t_0$, in order that
property \eqref{equation:bilip} may be applied. We therefore impose a condition
on $\gamma_1$ and $\gamma_2$, namely that
$$d_{t_0}(\gamma_1 \intersect P_{t_0}, 
          \gamma_2 \intersect P_{t_0})
   \ge \frac{C_6}{A}
$$
which implies, for all $t \le t_0$, that 
$$d_t(\gamma_1 \intersect P_t, \gamma_2 \intersect P_t) \ge C_6
$$
and so
\begin{align*}
\frac{1}{2K} \cdot d_t(\gamma_1 \cap P_t, \gamma_2 \cap P_t)
   &\le d'_{h(t)}(\gamma'_1 \cap P'_{h(t)},\gamma'_2 \cap P'_{h(t)})
   \\
   &\le 2K \cdot d_t(\gamma_1 \cap P_t, \gamma_2 \cap P_t)
\end{align*}
which implies
\begin{align*}
\frac{A}{2K} \alpha^{-t+t_0} d_{t_0}(\gamma_1 \intersect P_{t_0}, 
                                    \gamma_2 \intersect P_{t_0}) 
     &\le d'_{h(t)}(\gamma'_1 \intersect P'_{h(t)}, 
                    \gamma'_2 \intersect P'_{h(t)}) \\
     &\le 2BK \alpha^{-t+t_0} d_{t_0}(\gamma_1 \intersect P_{t_0}, 
                                     \gamma_2 \intersect P_{t_0})
\end{align*}

Next, applying the Exponential Lower Bound of Proposition
\ref{prop:growth} in $G_N$, for each $s \le s_0$ we have
$$ d'_{s}(\gamma'_1 \intersect P'_s,\gamma'_2 \intersect P'_s) \ge
  A \cdot \beta^{-s+s_0} d'_{s_0}(\gamma'_1 \intersect P'_{s_0},
                                  \gamma'_2 \intersect P'_{s_0})
$$
Taking $s=h(t)$, and using the fact that $s_0=h(t_0)$, this implies 
$$\beta^{-h(t)+h(t_0)} d'_{h(t_0)}(\gamma'_1 \intersect P'_{h(t_0)},
                                   \gamma'_2 \intersect P'_{h(t_0)}) 
   \le \frac{2BK}{A} \cdot \alpha^{-t+t_0} \cdot 
                       d_{t_0}(\gamma_1 \intersect P_{t_0},
                               \gamma_2 \intersect P_{t_0})
$$
Therefore,
$$\beta^{-h(t)+h(t_0)} d_{t_0}(\gamma_1 \intersect P_{t_0},
                               \gamma_2 \intersect P_{t_0}) 
\le \frac{4BK^2}{A} \alpha^{-t+t_0} d_{t_0}(\gamma_1 \intersect P_{t_0},
                                \gamma_2 \intersect P_{t_0})
$$
Now divide both sides by $d_{t_0}(\gamma_1\intersect P_{t_0}, 
\gamma_2 \intersect P_{t_0})$, and take logarithms,
obtaining
$$(-h(t) + h(t_0)) \log(\beta) \le \log\left(\frac{4BK^2}{A}\right) +
(-t+t_0) \log(\alpha)
$$
and so
$$h(t) \ge \frac{\log(\alpha)}{\log(\beta)} (t-t_0) + h(t_0) - 
\frac{\log\left(\frac{4BK^2}{A}\right)}{\log(\beta)}
$$
proving Claim \ref{claim:growth} and therefore completing the proof of
Proposition \ref{proposition:height:rigidity}.
\end{proof}

\subsection{Interlude: The induced boundary map}
\label{section:induced:boundary}

The \emph{upper boundary} $\bdy^u G_M$ is defined to be the leaf space
of the weak stable foliation; this leaf space is identified with $V^+$. The
\emph{lower boundary} $\bdy_\ell G_M$ is the leaf space of the weak
unstable foliation, identified with $V^-$. The {\em internal boundary}
$\bdy _\ext G_M$ is defined as 
$$
\bdy _\ext G_M = \bdy_\ell G_M \cross \bdy^u G_M = V^- \cross V^+ \approx 
\R^n / V^0
$$
which is identified with the leaf space of the center foliation. 

As a consequence of Proposition
\ref{proposition:foliations:preserved}, a quasi-isometry $\phi \from
G_M \to G_L$ which respects the transversely oriented horizontal
foliations induces a bijection
$$\bdy_\ext \phi \from \bdy_\ext G_M \to \bdy_\ext G_L
$$
which preserves the factors, that is,
$$\bdy_\ext\phi = \bdy_l \phi \cross \bdy^u \phi \from \bdy_\ell G_M 
\cross \bdy^u G_M \to \bdy_\ell G_L \cross \bdy^u G_L
$$

Recall the \nb{1}parameter family of metrics $d_{M,t}$ on $\R^n$ given
by the quadratic form $Q_{M,t} = (M^{-t})^T M^{-t}$. The internal
boundary $\bdy_\ext G_M$ is identified with $\R^n / V^0$ and with 
$V^- \cross V^+$, and we consider two \nb{1}parameter families of
metrics. 

First, regarding points of $\R^n / V^0$ as affine subspaces
parallel to $V^0$, there is a \nb{1}parameter family of Hausdorff
metrics induced from $d_{M,t}$ which we denote $dh_{M,t}$. Second,
restrict the action of $M^t$ to the subspace $V^- \cross V^+$ to get a
\nb{1}parameter subgroup of $\GL(V^- \cross V^+)$, and by choosing a
basis for $V^- \cross V^+$ we obtain a \nb{1}parameter subgroup
$\hat M^t$ of $\GL(k,\R)$, where $k$ is the dimension of 
$V^- \cross V^+$. We obtain a \nb{1}parameter family of metrics $d_{\hat 
M,t}$.
There is a canonical identification $V^- \cross V^+ \approx \R^n /
V^0$, and with respect to this identification the metrics $d_{\hat
M,t}$ and $dh_{M,t}$ are bilipschitz equivalent, with a uniform
bilipschitz constant independent of $t$.

Note that the absolute Jordan form of $\hat M$ is identical with the
nonunipotent part of the absolute Jordan form of $M$, and similarly
for $N$. 

\begin{lemma}
\label{lemma:parameterizedqi}
Given two \nb{1}parameter subgroups $M^t, N^t$ of $\GL(n,\R)$, for all
$K \ge 1$, $C,A \ge 0$, there exist $K' \ge 1, C' \ge 0$ with the
following properties. If $\phi \from G_M \to G_N$ is a $K,C$
quasi-isometry which coarsely respects the transversely oriented
horizontal foliations, with Hausdorff constant $A$, then for every $t
\in\R$ the induced map $\bdy_\ext\phi \from \bdy_\ext G_M \to \bdy_\ext
G_N$ is a $K',C'$ quasi-isometry from the metric $d_{\hat M,t}$ to the 
metric
$d_{\hat N,h(t)}$.
\end{lemma}

\begin{proof} With what we know, the proof is mostly a matter of chasing
through definitions.

The quasi-isometry $\phi$ is a bounded distance from a quasi-isometry
$\psi \from G_M \to G_N$ which takes the horizontal leaf $P_t$ to the horizontal leaf
$P'_{h(t)}$, and which simultaneously takes center leaves of $G_M$ to center leaves of
$G_N$. Now restrict the center foliations of $G_M, G_N$ to $P_t$,
$P'_{h(t)}$, and denote the respective leaf spaces as $Q_t, Q'_{h(t)}$. 

In order to apply Lemma \ref{lemma:path:subspace}, consider each horizontal leaf
$P_t$ of $G_M$ as a geodesic metric space, with respect to the Riemannian metric
induced by restriction from $G_M$. The inclusion map $P_t \inject G_M$ is evidently
$(1,0)$ coarsely lipschitz, and it is uniformly proper, with a uniformity function
$s(r) = a^r$ where $a>1$ is larger than the maximum of the absolute values of all
eigenvalues of $M$ and their multiplicative inverses. Note in particular that the
coarse Lipschitz constants and the uniformity functions of the maps $P_t \inject G_M$
depend only on $K,C,A$ and on the matrix $M$, but not on $t$. Similar remarks apply
to the inclusion map $P'_{h(t)} \inject G_N$. Applying Lemma
\ref{lemma:path:subspace}, restricting $\psi$ to $P_t$ results in a map $\psi_t
\from P_t \mapsto P'_{h(t)}$ which is a quasi-isometry. There is in turn an induced
map $\theta_t \from Q_t \mapsto Q'_{h(t)}$ which is a quasi-isometry with respect to
the associated Hausdorff metric. The quasi-isometry constants of the maps $\psi_t$
and $\theta_t$ depend only on $K,C,A$.

Now consider the coordinate identifications $G_M \approx \R^n \cross
\R$, $G_N\approx \R^n \cross \R$. By construction of the left invariant
metrics, for each $t$ the space $P_t$ is identified with $\R^n \cross
t \approx \R^n$ with the metric $d_{M,t}$, and the space $P'_{h(t)}$ is
identified with $\R^n$ with metric $d_{N,h(t)}$, and so the maps $\psi_t 
\from
\R^n \to \R^n$ are uniform quasi-isometries from $d_{M,t}$ to $d_{N,h(t)}$ 
for
all $t$. Also, $Q_t$ is identified with $\R^n / V^0_M$ with the
associated Hausdorff metric $dh_{M,t}$, and $Q'_{h(t)}$ is identified
with $\R^n / V^0_N$ with the associated Hausdorff metric $dh_{N,h(t)}$, 
and so the maps $\theta_t \from \R^n / V^0_M \to \R^n /
V^0_N$ are uniform quasi-isometries from $dh_{M,t}$ to $dh_{N,h(t)}$
for all $t$. This implies that $\theta_t \from V^-_M \cross V^+_M \to
V^-_N \cross V^+_N$ is a quasi-isometry from $d_{\hat M,t}$ to
$d_{\hat N,t}$ for all $t$. But for all $t$ the map $\theta_t$ is
identical to $\bdy_\ext\phi \from\bdy_\ext G_M \to \bdy_\ext G_N$,
proving the lemma.
\end{proof}

\subsection{Step 3: Reduction to Theorem \ref{theorem:parameter:rigidity}
on 1-parameter subgroup rigidity}
\label{section:reduction}

Assume the hypotheses of Theorem \ref{theorem:levelset:preserving},
namely that we have \nb{1}parameter subgroups $M^t, N^t$, and a
quasi-isometry $\phi \from G_M \to G_N$ which coarsely respects the
transversely oriented horizontal foliations. Applying Proposition
\ref{proposition:height:rigidity}, there is an induced time change of
the form $h(t) = mt+b$ with $m>0$. Applying Proposition
\ref{proposition:choices}, there is a horizontal respecting 
quasi-isometry $G_N \to G_{N^{m}}$ with an induced time change of
the form $s \mapsto s/m$. By composition we obtain a horizontal
respecting quasi-isometry $G_M \to G_{N^m}$ with an induced time change
of the form $t \mapsto t + b'$. Changing the coordinates in $G_M$ by a
translation of the time coordinate $t$, we have a horizontal
respecting quasi-isometry $G_M\to G_{N^m}$ for which the identity map
$t \mapsto t$ is an induced time change. Applying Lemma
\ref{lemma:parameterizedqi}, we obtain a map $\bdy_\ext \phi \from
\R^n \to \R^n$ which, for each $t$, is a $(K',C')$-quasi-isometry from $d_{\hat
M,t}$ to $d_{\hat N^m,t}$.

Now apply the following theorem (with $N$ in place of $N^m$), which will 
be proved in the next section:

\begin{theorem}[1-parameter subgroup rigidity] 
\label{theorem:parameter:rigidity}\quad\break
Let $M^t,N^t$ be \nb{1}parameter subgroups of $\GL(n,\R)$, such that
$M=M^1$ and $N=N^1$ have no eigenvalues on the unit circle. If 
there exists a bijection $f \from \R^n\to \R^n$ and constants $K \ge
1,C\ge 0$ such that for each $t\in \R$ and $p,q \in \R^n$ we have
$$-C+\frac{1}{K}d_{M,t}\bigl(p,q\bigr)\leq d_{N,t}(f(p),f(q))\leq K
d_{M,t}\bigl(p,q\bigr)+C
$$
then $M$ and $N$ have the same absolute Jordan form.  
\end{theorem}

Returning to the previous discussion, this theorem allows us to
conclude that $\hat M$ and $\hat N^m$ have the same absolute
Jordan form, and so the nonunipotent parts of the absolute Jordan
forms of $M,N^m$ are identical. We have already proved in
Corollary \ref{corollary:unipotent} that the unipotent parts are
identical, and so $M$ and $N^m$ have the same absolute Jordan
forms, finishing the proof of Theorem
\ref{theorem:horizontal:respecting}.
\qed


\section{Dynamics of $G_M$, Part II:\\ 1-parameter subgroup rigidity}
\label{section:oneparameter}

In this section we give a proof of Theorem
\ref{theorem:parameter:rigidity}. 

Let $M^t$, $N^t$ be \nb{1}parameter subgroups of $\GL(n,\R)$ with
no eigenvalues on the unit circle. Let $M^t = \barM^t P^t$, $N^t =
\barN^t Q^t$ be the real Jordan forms, so $\barM$ and $\barN$ have
all positive eigenvalues, none equal to $1$. Let $f \from \R^n \to
\R^n$ be a bijection which satisfies
\begin{equation}\label{equation:qi}
-C+\frac{1}{K}d_{M,t}\bigl(p,q\bigr)\leq d_{N,t}(f(p),f(q))\leq K
d_{M,t}\bigl(p,q\bigr)+C
\end{equation}
for all $t \in \R, p,q\in \R^n$.

The bijection $f \from \R^n \to \R^n$ must in fact be a homeomorphism. To
see why, for each $p\in \R^n$, $R>0$, $T>0$ let
$$F_{p,R}(T)=\{q\in\R^n \suchthat d_{M,t}(p,q)<R \quad\text{for all}\quad
t\in (-T,T)\}
$$
In other words, $F_{p,R}(T)$ is the intersection of open balls of
radius $R$ about $p$ in each of the metrics $d_{M,t}$, for $t \in
(-T,T)$.  Since the eigenvalues of $\barM$ are all positive real
numbers, none equal to $1$, it follows from Proposition
\ref{prop:growth} that for each $p \in \R^n$ and each $R>0$ the
collection of sets $F_{p,R}(T)$ as $T$ ranges in $(0,\infty)$ is a
neighborhood basis for $p$, in the standard topology on $\R^n$.  We
define a similar neighborhood basis using matrix $N$, denoted
$G_{p,R}(T)$.  Since $f(F_{p,R}(T))\subseteq G_{f(p),KR+C}(T)$ for
each $p\in\R^n,R>0,T>0$, it follows that $f$ is continuous.  The same
argument applies to $f^{-1}$, and so $f$ is a homeomorphism.

The idea of the proof of Theorem \ref{theorem:parameter:rigidity} is
to show that $f$ respects certain ``flags of foliations'' which are
closely related to the Jordan decompositions of $\R^n$ with respect to
$M^t$ and $N^t$. We begin by setting up the notation needed to define
and study these foliations.

\begin{Definition}[Flags of foliations]
If $V$ is a vector subspace of $\R^n$, define a foliation $\F(V)$ of
$\R^n$ whose leaves are the affine subspaces of $\R^n$ parallel to
$V$. Given a flag of subspaces $V_1 \subset \cdots \subset V_r$, it
follows that if $1 \le i < j \le r$ then each leaf of $\F(V_i)$ is
contained in some leaf of $\F(V_j)$; we denote this relation by saying
that $\F(V_1) \prec \cdots \prec \F(V_r)$ is a \emph{flag of
foliations} of $\R^n$. 
\end{Definition}

Recall the root space decompositions of $\R^n$ with respect to $\barM$ and
$\barN$. We denote the eigenvalues of $\barM$ and $\barN$ by 
$$0<\mu^-_m< \cdots <\mu^-_1<1<\mu^+_1< \cdots <\mu^+_r$$ 
and 
$$0<\nu^-_n< \cdots <\nu^-_1<1<\nu^+_1<\cdots <\nu^+_s$$ 
respectively. The corresponding root space decompositions are
denoted
$$V^-_m \oplus \cdots \oplus V^-_1 \oplus V^+_1 \oplus \cdots
\oplus V^+_r
$$
and
$$W^-_n \oplus \cdots \oplus W^-_1 \oplus W^+_1 \oplus \cdots
\oplus W^+_s.
$$
As in section \ref{section:liegroup} we set
\begin{align*}
V^- = V^-_m \oplus \cdots \oplus V^-_1, &\qquad 
   V^+ = V^+_1 \oplus \cdots \oplus V^+_r \\
W^- = W^-_n \oplus \cdots \oplus W^-_1, &\qquad
   W^+ = W^+_1 \oplus \cdots \oplus W^+_s
\end{align*}
Define the \emph{root space flags}
\begin{align*}
U^-_i &= V^-_i \oplus \cdots \oplus V^-_1, \quad i =  1,\ldots,m \\
U^+_j &= V^+_1 \oplus \cdots \oplus V^+_j, \quad j =  1,\ldots,r \\
Y^-_i &= W^-_i \oplus \cdots \oplus W^-_1, \quad i =  1,\ldots,n \\
Y^+_j &= W^+_1 \oplus \cdots \oplus W^+_j, \quad j =  1,\ldots,s
\end{align*} 
and by convention we take $U^-_0$, $U^+_0$, $Y^-_0$, $Y^+_0$ each to be the
trivial subspace.  Associated to the root space flags we have
\emph{root space foliation flags}
\begin{align*}
\F(U^-_1) &\prec \cdots \prec \F(U^-_m) = \F(V^-) \\
\F(U^+_1) &\prec \cdots \prec \F(U^+_r) = \F(V^+) \\
\F(Y^-_1) &\prec \cdots \prec \F(Y^-_n) = \F(W^-) \\
\F(Y^+_1) &\prec \cdots \prec \F(Y^+_s) = \F(W^+)
\end{align*}

\subsubsection*{Step 1: $f$ respects contracting and expanding
foliations.} First we show that $f(\F(V^-)) = \F(W^-)$ and $f(\F(V^+))
= \F(W^+)$. 

Given $p,q \in \R^n$ we have the following chain of equivalences:
\begin{enumerate}
\item $p,q$ are in the same leaf of $\F(V^+)$.
\item $d_{M,t}(p,q) = \| M^{-t}(p-q) \| \to 0$ as $t \to
+\infinity$.
\item $d_{M,t}(p,q)$ is bounded for $t \in [0,+\infinity)$.
\item $d_{N,t}(f(p),f(q))$ is bounded for $t \in [0,+\infinity)$.
\item $d_{N,t}(f(p),f(q)) = \| N^{-t}(f(p)-f(q)) \| \to 0$ as $t \to
+\infinity$. 
\item $f(p), f(q)$ are in the same leaf of $\F(W^+)$.
\end{enumerate}
The equivalence of (1--3) follows from Proposition
\ref{prop:growth}, and similarly for (4--6). The equivalence of
(3) and (4) follows from \eqref{equation:qi}. This shows $f(\F(V^+)) = \F(W^+)$. A
similar argument with $t \in (-\infinity,0]$ shows $f(\F(V^-)) =
\F(W^-)$.

\subsubsection*{Step 2: $f$ respects root space foliation flags.} 
Next we show:

\begin{claim}
\label{claim:rootflags}
$f \from \R^n \to \R^n$ respects the root space foliation flags, and
corresponding root spaces have the same eigenvalues. More precisely we have:
\begin{enumerate}
\item $r=s$.
\item $\mu^+_j = \nu^+_j$ for $j=1,\ldots,r$.
\item $f(\F(U^+_j)) = \F(Y^+_j)$ for $j=1,\ldots,r$.
\item $m=n$.
\item $\mu^-_i = \nu^-_i$ for $i=1,\ldots,m$.
\item $f(\F(U^-_i)) = \F(Y^-_i)$ for $i=1,\ldots,m$.
\end{enumerate}
It follows that $M,N$ have the same eigenvalues with the same multiplicities.
\end{claim}

We give the proof of 1,2,3; the proof of 4,5,6 is similar.

We know by Step 1 that $f(\F(V^+)) = \F(W^+)$. Consider points $p,q$
in the same leaf of $\F(V^+)$, so $f(p),f(q)$ are in the same leaf of
$\F(W^+)$. From Proposition \ref{prop:growth} it follows that
as $t \to -\infinity$ both of the quantities $d_{M,t}(p,q)$ and
$d_{N,t}(f(p),f(q))$ approach $+\infinity$. It follows that for
sufficiently large $t$, in the inequality \eqref{equation:qi} we can absorb the
additive constant $C$, yielding
\begin{equation}\label{equation:coarsebilip}
\frac{1}{K+1} d_{M,t}(p,q) \le d_{N,t}(f(p),f(q)) \le (K+1) d_{M,t}(p,q)
\end{equation}
Define displacement vectors $v=p-q$, $w=f(p)-f(q)$. Taking natural
logarithms, dividing by $t$, and taking limsup, we have
\begin{align}
\limsup_{t \to -\infinity} \frac{\log\bigl(d_{M,t}(p,q)\bigr)}{t} 
  &= \limsup_{t \to -\infinity} \frac{\log\bigl(d_{N,t}(f(p),f(q))\bigr)}{t} \notag\\
\label{equation:equalrates}
\limsup_{t \to +\infinity} \frac{\log \| M^t v \|}{t}
  &= \limsup_{t \to +\infinity} \frac{\log \| N^t w \|}{t}
\end{align}
To evaluate these limits, let $I(p,q)=I(v)$ be the unique integer such
that 
$$v \in U^+_{I(v)} - U^+_{I(v)-1}
$$
or, equivalently, the unique integer such that $p,q$ are in the same
leaf of $\F(U^+_{I(p,q)})$ but not in the same leaf of
$\F(U^+_{I(p,q)-1})$ (recall the convention that $U^+_0 = 0$, and so
$I(p,q)=0$ if and only if $p=q$). Define $J(f(p),f(q)) = J(w)$ similarly by
$$w \in Y^+_{J(w)} - Y^+_{J(w)-1}
$$
Applying Proposition \ref{prop:growth} we have
\begin{align*}
\limsup_{t \to +\infinity} \frac{\log \| M^t v \|}{t} &= \mu^+_{I(v)} \\
\limsup_{t \to +\infinity} \frac{\log \| N^t w \|}{t} &= \nu^+_{J(w)}
\end{align*}
and so by \eqref{equation:equalrates} we have
$$
\mu^+_{I(p,q)} = \mu^+_{I(v)} = \nu^+_{J(w)} = \nu^+_{J(f(p),f(q))}
$$
Since $f$ is a bijection from each leaf of $\F(V^+)$ to some leaf of
$\F(W^+)$, items (1) and (2) of Claim \ref{claim:rootflags} now
follow, and it also follows that
$$I(p,q) = J(f(p),f(q))
$$
for all $p,q$ contained in the same leaf of $\F(V^+)$.

We now prove item (3) of Claim \ref{claim:rootflags} by induction on
$j$. If $p,q$ are in the same leaf of $\F(U^+_1)$ then $I(p,q)=1$ and
so $J(p,q)=1$ which implies that $f(p),f(q)$ are in the same leaf of
$\F(Y^+_1)$. A similar argument with $f^\inverse$ proves that
$f(\F(U^+_1)) = \F(Y^+_1)$, proving the base step of the
induction. Now assume that $f(\F(U^+_j)) = \F(Y^+_j)$, and suppose
$p,q$ are in the same leaf of $\F(U^+_{j+1})$. There are two cases to
consider. If $p,q$ lie in the same leaf of $\F(U^+_j)$ then by the
induction hypothesis $f(p),f(q)$ lie in the same leaf of $\F(Y^+_j)$,
in particular they lie in the same leaf of $\F(Y^+_{j+1})$. If $p,q$
do not lie in the same leaf of $\F(U^+_j)$ then $I(p,q) = j+1$ and so
$J(f(p),f(q)) = j+1$ and so $f(p),f(q)$ lie on the same leaf of
$\F(Y^+_{j+1})$. A similar argument with $f^\inverse$ shows that
$f(\F(U^+_{j+1}))=Y^+_{j+1}$, completing the induction.

As mentioned earlier, (4--6) are proved similarly, completing the
proof of Claim \ref{claim:rootflags}.

\subsubsection*{Step 3: $f$ respects Jordan foliation flags.}

From Step 2, for each fixed $j=1,\ldots,r$ the matrices $M,N$ have
$\mu^+_j$ root spaces $V^+_j, W^+_j$ respectively. As part of their root
space flags we have
\begin{align*}
U^+_j &= U^+_{j-1} \oplus V^+_j \\
Y^+_j &= Y^+_{j-1} \oplus W^+_j
\end{align*}
Let $c_j$ be the index of nilpotency of $\mu_j \cdot I - M$, and let
$d_j$ be the index of nilpotency of $\mu_j \cdot I - N$. Then we have
Jordan filtrations
\begin{align*}
V^+_{j,0} \subset \cdots \subset V^+_{j,c_j} &= V^+_j \\
W^+_{j,0} \subset \cdots \subset W^+_{j,d_j} &= W^+_j
\end{align*}
and we set $U^+_{j,k} = U^+_{j-1} \oplus V^+_{j,k}$ and $Y^+_{j,k} =
Y^+_{j-1} \oplus W^+_{j,k}$, yielding subspace flags
\begin{align*}
U^+_{j-1} \subset U^+_{j,0} \subset \cdots \subset U^+_{j,c_j-1} &=
   U^+_j \\
Y^+_{j-1} \subset Y^+_{j,0} \subset \cdots \subset Y^+_{j,d_j-1} &= Y^+_j
\end{align*}
Corresponding to these subspace flags are foliation flags, 
\begin{align*}
\F(U^+_{j-1}) \prec \F(U^+_{j,0}) \prec \cdots \prec \F(U^+_{j,c_j-1}) &= \F(U^+_j) \\
\F(Y^+_{j-1}) \prec \F(Y^+_{j,0}) \prec \cdots \prec \F(Y^+_{j,d_j-1}) &= \F(Y^+_j)
\end{align*}
called the \emph{expanding Jordan foliation flags} associated to
the corresponding root space foliations $\F(U^+_j)$,
$\F(Y^+_j)$ respectively. The \emph{contracting Jordan foliation
flags} associated to each root space foliation $\F(U^-_i)$,
$\F(Y^-_i)$ are similarly defined.

\begin{claim}
\label{claim:Jordanflags}
$f \from \R^n \to \R^n$ respects the Jordan foliation flags associated to corresponding root space foliations. More
precisely, for each $j=1,\ldots,r$ we have:
\begin{enumerate}
\item $c_j = d_j$.
\item $f(\F(U^+_{j,k})) = \F(Y^+_{j,k})$ for $k=0,\ldots,c_j-1$.
\end{enumerate} 
and similarly for the contracting Jordan foliation flags.
\end{claim}

From this claim, for each $j=1,\ldots,r$ it immediately follows
that $\barM, \barN$ have the same Jordan blocks with eigenvalue
$\mu^+_j$, and so the expanding parts of the Jordan forms for
$\barM, \barN$ are identical; similarly for the contracting parts.
Since $M,N$ have no eigenvalues on the unit circle, it now follows
that $M,N$ have the same absolute Jordan forms, completing the
proof of Theorem \ref{theorem:parameter:rigidity}.

\begin{proof}[Proof of Claim \ref{claim:Jordanflags}]
Consider $p,q \in \R^n$ in the same leaf of $\F(U^+_j)$ but not in the
same leaf of $\F(U^+_{j-1})$, and so $f(p),f(q)$ are in the same leaf
of $\F(Y^+_j)$ but not in the same leaf of $\F(Y^+_{j-1})$. Define
displacement vectors $v=p-q$, $w=f(p)-f(q)$, and so $v \in U^+_j -
U^+_{j-1}$ and $w \in Y^+_j - Y^+_{j-1}$. We know that
$$\limsup_{t \to +\infinity} \frac{\log \| M^t v \|}{t} = \frac{\log
\| N^t w \|}{t} = \mu^+_j
$$
We also know that \eqref{equation:coarsebilip} is true for $t$
sufficiently close to $-\infinity$, and so for $t$ sufficiently
close to $+\infinity$ we have
\begin{equation}\label{equation:samegrowth}
\frac{1}{K+1} \| M^{t} v \| \le \| N^{t} w \| 
                               \le (K+1) \| M^{t} v \| 
\end{equation}

By induction on $k=0,1,\ldots$, we shall prove that $v \in U^+_{j,k}$ if and
only if $w \in Y^+_{j,k}$, or equivalently that $f(\F(U^+_{j,k})) =
\F(Y^+_{j,k})$. 

For the basis step $k=0$, divide the inequality \eqref{equation:samegrowth} by $\mu^t$ to
obtain, for all $t$ sufficiently close to $+\infinity$:
\begin{equation}\label{equation:samegrowth2}
\frac{1}{K+1} \frac{\| M^{t} v \|}{\mu^t} \le \frac{\| N^{t} w \|}{\mu^t} 
                               \le (K+1) \frac{\| M^{t} v \|}{\mu^t} 
\end{equation}
By the Exponential Lower Bound and the Exponential$\cdot$Polynomial
Upper Bound of Proposition \ref{prop:growth}, the quantity
$\displaystyle\frac{\| M^t v \|}{\mu^t}$ is bounded for $t \ge 0$ if
and only if $v \in U^+_{j,0}$; and the quantity $\displaystyle\frac{\|
N^t w \|}{\mu^t}$ is bounded on $t \ge 0$ if and only if $w \in
Y^+_{j,0}$. However by \eqref{equation:samegrowth2} the boundedness of these two
quantities on $t \ge 0$ are equivalent.

For the induction step, assume that $f(\F(U^+_{j,k-1})) =
\F(Y^+_{j,k-1})$, that is, $v \in U^+_{j,k-1}$ if and only if $w \in
Y^+_{j,k-1}$. We must prove that $v \in U^+_{j,k} - U^+_{j,k-1}$ if
and only if $w \in Y^+_{j,k} - Y^+_{j,k-1}$. From 
\eqref{equation:samegrowth}, for $t$ sufficiently close to $+\infinity$ we have
\begin{gather}
\label{equation:degreek}
\frac{1}{K+1} \frac{\| M^{t} v \|}{\mu^t t^k} 
                 \le \frac{\| N^{t} w \|}{\mu^t t^k} 
                 \le (K+1) \frac{\| M^{t} v \|}{\mu^t t^k}\\
\intertext{and}
\label{equation:degreek-1}
\frac{1}{K+1} \frac{\| M^{t} v \|}{\mu^t t^{k-1}} 
                 \le \frac{\| N^{t} w \|}{\mu^t t^{k-1}} 
                 \le (K+1) \frac{\| M^{t} v \|}{\mu^t t^{k-1}} 
\end{gather}
By the Exponential$\cdot$Polynomial Upper and Lower Bounds of
Proposition \ref{prop:growth}, the following two statements are
equivalent:
\begin{itemize}
\item[(1)]
$v \in U^+_{j,k} - U^+_{j,k-1}$
\item[(2)] 
For $t\ge 0$, the quantity $\displaystyle\frac{\| M^{t} v \|}{\mu^t
t^k}$ is bounded, but the quantity $\displaystyle\frac{\| M^{t} v
\|}{\mu^t t^{k-1}}$ is not bounded.
\end{itemize}
Similarly, the following two statements are equivalent:
\begin{itemize}
\item[(3)] 
$w \in Y^+_{j,k} - Y^+_{j,k-1}$.
\item[(4)] 
For $t \ge 0$, the quantity $\displaystyle\frac{\| N^{t} w \|}{\mu^t
t^k}$ is bounded, but the quantity $\displaystyle\frac{\| N^{t} w
\|}{\mu^t t^{k-1}}$ is not bounded.
\end{itemize}
But by inequalities \eqref{equation:degreek} and
\eqref{equation:degreek-1}, statements (2) and (4) are equivalent,
and so statements (1) and (3) are equivalent, completing the
inductive proof of item 2 of Claim \ref{claim:Jordanflags} for all
$k \ge 0$.

The foliation flag $\F(U^+_{j,0}) \prec \cdots \prec \F(U^+_{j,k})
\prec \cdots$ must terminate at $\F(U^+_j)$ for the same value of $k$
for which the flag $\F(Y^+_{j,0}) \prec \cdots \prec \F(Y^+_{j,k})
\prec \cdots$ terminates at $\F(Y^+_j)$, proving that $c_j=d_j$, and
completing the proof of Claim \ref{claim:Jordanflags}.
\end{proof}

Our proof of Theorem \ref{theorem:parameter:rigidity} actually  provides
for some regularity of $f$. We record the statement here, although it is
not used at all in this paper.

\begin{proposition}[Regularity]
\label{proposition:regularity}
With the assumptions as in Theorem \ref{theorem:parameter:rigidity}, $f$
is a homeomorphism which respects the contracting and expanding root
space foliation flags of $\overline M, \overline N$, and for each
corresponding pair of root space foliations $f$ also respects the
associated Jordan foliation flags.\qed
\end{proposition}

\begin{Remark} Even stronger regularity properties should hold.
For instance, $f$ should satisfy lipschitz conditions in directions
parallel to a root space,  by arguments similar to the results of
\cite{FarbMosher:BSOne}. Understanding what happens transverse to root
spaces will require new ideas.
\end{Remark}


\section{Quasi-isometries of $\Gamma_M$ via Coarse Topology}
\label{section:coarsetop}

Recall the notation for abelian-by-cyclic Lie groups: given $M \in
\GL_\cross(m,\R)$, a \nb{1}parameter subgroup $M^t\subset \GL(m,\R)$ with
$M^1=M$ determines a Lie group denoted $G_M = \R^m \semidirect_M \R$.

This entire section will be devoted to a proof of the following. 

\begin{proposition}[Induced quasi-isometries of $G_M$]
\label{proposition:induced}
\label{proposition:vhpreserved}
Consider integral matrices $M \in \GL_\cross(m,\R)$, $N \in
\GL_\cross(n,\R)$ and suppose that $\det M, \det N > 1$. If there exists
a quasi-isometry $f\from\Gamma_M\to \Gamma_N$ then $m=n$ and there exists
a quasi-isometry $\phi\from G_M\to G_N$ which coarsely respects horizontal
foliations and their transverse orientations. Furthermore, all associated
constants for $\phi$ depend only on those for~$f$.
\end{proposition}

\subsection{A geometric model for $\Gamma_M$}
\label{section:model}

Let $M \in \GL_\cross(m,\R)$ be an integral matrix lying on a
\nb{1}parameter subgroup $M^t$ of $\GL(m,\R)$ with $M^1=M$ and with
associated Lie group $G_M$. We assume that $\det M > 1$ and we denote
$d=\det M$.
  
We start by constructing a contractible, $(m+1)$-dimensional
metric complex $X_M$ on which $\Gamma_M$ acts properly
discontinuously and cocompactly by isometries, and so the group
$\Gamma_M$ will be quasi-isometric to the geodesic metric space
$X_M$.

The description of $\Gamma_M$ as an ascending HNN extension shows
that $\Gamma_M$ is the fundamental group of the mapping torus of
an injective endomorphism of the $m$-dimensional torus.  Let $X_M$
be the universal cover of this mapping torus.  Topologically,
there is a fibration
$$\xymatrix{
\R^{n-1} \ar[r] & X_M \ar[d] \\
                & T_M \\
}
$$
where $T_M$ is the homogeneous directed tree with one edge coming into
each vertex and $d=\det M$ edges going out of each vertex. Hence $X_M$ 
is a topological product $X_M\approx \R^{n-1}\times T_M$.

The action of 
$\Gamma_M$ on $X_M$ by deck tranformations induces an action of 
$\Gamma_M$ on $T_M$. This action is equivalent to the 
usual action of the HNN extension 
$\Gamma_M$ on its Bass-Serre tree $T_M$.

Before constructing a metric on $X_M$, let us describe the essential 
properties of such a metric. These are best described by 
giving the isometry types of natural subcomplexes of $X_M$.

\begin{Definition}[Doubled horoballs]
We define a {\em doubled $G_M$ horoball}, denoted $H_M$, to be the
metric space obtained by identifying two copies of $\{(x,t)\in G_M
\suchthat t\geq 0\}$ along $\{(x,0)\in G_M\}$, endowed with the
path metric. 
\end{Definition}

\begin{Definition}[Hyperplanes in $X_M$]
Let $P_{\ell}=\pi_M^{-1}(\ell)$, where $\ell$ is a bi-infinite line in
the directed tree $T_M$.  We call $P_\ell$ a {\em hyperplane} in $X_M$. 
There are two
cases to consider:
\begin{itemize}
\item $\ell$ is coherently oriented in $T_M$.  In this case 
$P_\ell$ is isometric to $G_M$, and we call $P_\ell$ a {\em
coherent hyperplane} in $X_M$.

\item $\ell$ is not coherently oriented in $T_M$, and thus 
switches orientation precisely once.  In this case $P_\ell$ is
isometric to $H_M$, and we call $P_\ell$ an {\em incoherent hyperplane} 
in $X_M$.
\end{itemize}
\end{Definition}

This definition nearly determines a metric on $X_M$.  To specify a
metric on $X_M$, one proceeds as follows.  Fix a path metric on
$T_M$ so that each edge has length $1$.  Fix a base vertex on
$T_M$. These choices determine a unique height function $T_M
\to\R$ taking the base vertex to the origin, and taking each edge
to a segment of length $1$ via an orientation preserving isometry.
We have also defined a height function $G_M \to \R$.  Note that the
height function on $G_M$ was previously called the ``time function''; 
we will use both terms.

The complex $X_M$ is the fiber product of the two height functions
$T_M \to \R$, $G_M \to \R$, as shown in the following diagram:
$$\xymatrix{
& X_M \ar[dl]_{g_M} \ar[dd] \ar[dr]^{\pi_M} \\
{G_M} \ar[dr] & & T_M \ar[dl] \\
& {\reals}
}
$$  
There are induced projections $g_M \from X_M \to G_M$ and $\pi_M
\from X_M \to T_M$, and an induced height function $X_M \to \R$.
There is a unique path metric on $X_M$ so that each continuous
cross section $G_M \to X_M$ of $g_M$ is a path-isometric
embedding; and hence each coherent hyperplane in $X_M$ is an
isometrically embedded copy of $G_M$.

\begin{Definition}[horizontal leaf]
A {\em horizontal leaf} $L$ in $X_M$ is a subset of the form 
$L=\pi_M^{-1}(v)$ where $v\in T_M$.  
\end{Definition}

Note that the collection of horizontal leaves on $X_M$, equipped with the 
Hausdorff metric, forms a metric space which is isometric to $T_M$ via the 
projection map $\pi_M \from X_M \to T_M$. 

Note that each hyperplane in $X_M$ comes equipped with a foliation
by horizontal leaves.  For coherent hyperplanes $P$ in $X_M$,
which are isometric to $G_M$, the notion of horizontal leaf in
$P$ coincides with that of a horizontal leaf in $G_M$, given in \S
\ref{section:horizontaldefs}.

\subsection{Proof of Proposition \ref{proposition:induced} on induced
quasi-isometries of $G_M$}
\label{section:coherent}

Let $M,N$ be as in the statement of the proposition. 

We begin by showing that $M$ and $N$ have the same size. Suppose that $M
\in\GL(m,\R)$ and $N \in \GL(n,\R)$. In \S\ref{section:model} we
constructed finite classifying spaces for $\Gamma_M$ and $\Gamma_N$ of
dimensions $m+1, n+1$ respectively, and by Lemma 5.2 of
\cite{FarbMosher:BSTwo} these numbers are the virtual cohomological
dimensions of $\Gamma_M, \Gamma_N$. By a result of Block-Weinberger
\cite{BlockWeinberger} and Gersten \cite{Gersten:dimension}, virtual
cohomological dimension is a quasi-isometry invariant for groups with
finite classifying spaces. It follows that $m=n$.

Now $\Gamma_M$ acts properly discontinuously, freely, and cocompactly on
$X_M$. This action is by isometries, because $\Gamma_M$ acts on $G_M$, on
$T_M$, and on $\R$ by isometries, and the fiber product diagram is
equivariant with respect to these actions. It follows that $\Gamma_M$ in
any word-metric is quasi-isometric to $X_M$. Henceforth we will freely
interchange $\Gamma_M$ and $X_M$ when discussing quasi-isometry type. The
same discussion applies to $\Gamma_N$ and $X_N$, and so the
quasi-isometry $f\from \Gamma_M \to \Gamma_N$ gives a quasi-isometry
(perhaps with bigger constants) $f \from X_M\to X_N$.

Proposition \ref{proposition:induced} generalizes the case when
$M$ and $N$ are $1\times 1$ matrices, done in \S 4 and \S 5 of
\cite{FarbMosher:BSOne}. The proof here is more difficult, and the steps
must be proved in different order. In Steps 1 and 2 we prove (in a more
general context; see Theorem \ref{theorem:horizontal}) that a
quasi-isometry $X_M \to X_N$ coarsely respects hyperplanes and horizontal
sets. However, we must still distinguish between coherent and incoherent
hyperplanes. This is easy in the $1
\cross 1$ case handled in \cite{FarbMosher:BSOne}, where $G_M$ and $G_N$
are (scaled versions of) $\hyp^2$, and a doubled $\hyp^2$ horoball is
evidently not quasi-isometric to $\hyp^2$. In general we are unable to
distinguish the quasi-isometry types of coherent and incoherent
hyperplanes. To get around this, in Step 3, Proposition
\ref{proposition:area}, we prove that there is no horizontal respecting
quasi-isometry between a coherent and an incoherent hyperplane.

\paragraph*{Step 1. Quasi-isometrically embedded hyperplanes are
close to hyperplanes:} Given integral matrices $M,N \in
\GL_\cross(n,\R)$, if $P=G_M$ or $H_M$, then for all $K \ge 1, C \ge 0$
there exists $A \ge 0$ such that if $\phi \from P \to X_N$ is a
$K,C$-quasi-isometric embedding then there is a unique hyperplane $Q
\subset X_N$ with $d_\Haus(\phi(P),Q)\leq A$.

\bigskip

This was proved for $1 \cross 1$ matrices in \cite{FarbMosher:BSOne}. Our
proof of Step 1, while following the same outline as the $1 \cross 1$
case, will actually apply in a much broader setting. The generalized
versions of Steps 1 and 2, given in Theorem
\ref{theorem:metricfibration} and Theorem \ref{theorem:horizontal}, are
used for example in \cite{FarbMosher:sbf} to study surface-by-free groups,
and also in \cite{MosherSageevWhyte} to prove quasi-isometric rigidity
theorems for various ``homogeneous'' graphs of groups (see the remark
after Theorem~\ref{theorem:horizontal}). 

The generalization of Step 1 given in Theorem
\ref{theorem:metricfibration} will require moving from the category of
quasi-isometric embeddings into the category of uniformly proper
embeddings. After a fair amount of work to establish the new setting, we
then quote some theorems of coarse algebraic topology and follow the
proof of \cite{FarbMosher:BSOne}.

Consider a finite graph $\Gamma$ of finitely generated groups;
each edge $e$ is oriented, with initial and final vertices $i(e)$,
$f(e)$. We say that $\Gamma$ is \emph{geometrically homogeneous} if each
edge-to-vertex injection is a quasi-isometry with respect to the word
metric space, or equivalently, has finite index image. Ideally we would
like to have a version of Step 1 for any geometrically homogeneous graph of
groups in which each vertex and edge group is the fundamental group of a
closed, aspherical $n$-manifold, or even more generally, an
$n$-dimensional \Poincare\ duality group. This should come from a more
careful reading of results in coarse algebraic topology such as
\cite{KapovichKleiner:duality}, but meanwhile we will use Theorems
\ref{theorem:coarse:separation} and \ref{theorem:packing}, which require
us to impose additional assumptions on $\Gamma$.

Suppose that we have a category $\mathcal C$ of aspherical, closed, smooth
manifolds such that $\mathcal C$ is closed under finite coverings and
satisfies \emph{smooth rigidity}, meaning that any homotopy equivalence
between manifolds in $\mathcal C$ is homotopic to a diffeomorphism. Such
categories include: the $n$-torus, $n \ge 1$; hyperbolic surfaces; all
other irreducible, nonpositively curved, locally symmetric spaces, by
Mostow's Rigidity Theorem \cite{Mostow:rigidity}; solvmanifolds, by
earlier work of Mostow \cite{Mostow:solvmanifolds}; nilmanifolds, by
still earlier work of Malcev \cite{Malcev:nilmanifolds}; and various
generalizations due to Farrell and Jones \cite{FarrellJones:rigidity},
\cite{FarrellJones:infrasolv}.

We shall assume that $\Gamma$ is a geometrically homogeneous graph of
groups where each vertex group $\Gamma_v$ is the fundamental group of a
manifold $M_v$ in the category $\mathcal C$. Construct a \emph{graph of
aspherical manifolds} $M_\Gamma$, with fundamental group $\pi_1
\Gamma$, as follows. For each edge $e$, the two injections $\Gamma_e \to
\Gamma_{i(e)}$, $\Gamma_e \to \Gamma_{t(e)}$ determine two finite covering
spaces of $M_v$ each of whose fundamental group is identified with
$\Gamma_e$, and so we obtain a diffeomorphism between the two covering
spaces; identify these covering spaces and let $M_e$ be the resulting
smooth manifold. We have smooth, finite covering maps $M_e \to M_{i(e)}$,
$M_e \to M_{t(e)}$ inducing the corresponding edge-to-vertex group
injections. Form $M_\Gamma$ from the disjoint union
$$\left(\Union_v M_v\right)\union \left(\Union_e M_e\cross e\right)
$$
by gluing $M_e \cross i(e)$ to $M_{i(e)}$ and $M_e \cross f(e)$ to
$M_{f(e)}$ via the finite covering maps $M_e \to M_{i(e)}$ and $M_e \to
M_{f(e)}$. From the construction of $M_\Gamma$ we obtain a map $M_\Gamma
\to\Gamma$, such that each fiber $M_x$, $x\in \Gamma$, is a manifold in
the category~$\mathcal C$.

Let $X_\Gamma$ be the universal cover of $M_\Gamma$. There is a
$\Gamma$-equivariant fiber bundle $X_\Gamma \to T_\Gamma$ over the
Bass-Serre tree $T_\Gamma$ of $\Gamma$, whose fiber is a contractible
$n$-manifold. Any geodesic metric on $M_\Gamma$ lifts to a
$\pi_1\Gamma$-equivariant geodesic metric on $X_\Gamma$. Smoothness
allows us to impose additional geometric structure on $X_\Gamma$ which we
now describe.

A geodesic metric space is \emph{proper} if closed balls are compact. A
\emph{bounded-geometry metric simplicial complex} is a simplicial complex
$\Sigma$ equipped with a proper, geodesic metric such that for some
constants $0<C_1<C_2$ each positive dimensional simplex has diameter
between $C_1$ and $C_2$, and for some constant $C > 0$ the link of each
simplex has $\le C$ simplices. A subset $S$ of $\Sigma$ is
\emph{rectifiable} if for any $p,q \in S$ there exists a path in $S$
between $p$ and $q$ which is rectifiable in $\Sigma$, and which has the
shortest $\Sigma$-length among all paths in $S$ between $p$ and
$q$. The length of such a path defines a geodesic metric on $S$. A
\emph{$D$-homotopy} in $\Sigma$ is a homotopy whose tracks all have
diameter $\le D$. The space $\Sigma$ is \emph{uniformly contractible} if
there exists a function $\delta\from [0,\infinity)\to [0,\infinity)$,
such that for every bounded subset $S\subset \Sigma$, the inclusion map $S
\inject \Sigma$ is \nb{$\delta\bigl(\diam(S)\bigr)$}homotopic to a
constant map. More precisely we say that $\Sigma$ is
\emph{$\delta$-uniformly contractible}.

Let $T$ be a bounded geometry, metric simplicial tree, let $X$ be a
proper, geodesic metric, and let $\pi \from X\to T$ be a surjective map.
Denote $X_A = \pi^\inv(A)$ for each $A\subset T$. The map $\pi$ is called
a \emph{metric fibration} if:

\begin{enumerate}
\item $X$ is a uniformly contractible, bounded geometry, metric
simplicial complex.
\label{ItemUnifContr}
\item For each subtree $T' \in T$, the subset $X_{T'}$ is a subcomplex of
$X$ and is rectifiable in $X$.
\label{ItemPathIsometry}
\item For each $t \in T$ the subspace $X_t$ is uniformly
contractible and is a bounded geometry, metric simplicial complex, with
bounded geometry constants and uniform contractiblity data independent of
$t$.
\label{ItemHorizUnifContr}
\item The map $\pi \from X \to T$ is distance nonincreasing.
\label{ItemNonincreasing}
\item There is a homeomorphism $\Theta\from X \to F \cross T$ such that:
\begin{enumerate}
\item For all $t \in T$, $\Theta(X_{t})=F\cross t$.
\item For all $x \in F$, the map $T \to x \cross T
\xrightarrow{\Theta^\inv} X$ is a locally isometric embedding.
\label{ItemTreePathIsometry}
\item There exists $K \ge 1$ such that for all edges $e$ of $T$ and $t
\in e$, the retraction $r \from e \to t$ induces a projection $X_e
\xrightarrow{\Theta} F \cross e \xrightarrow{\Id \cross r} F \cross t
\xrightarrow{\Theta^\inv} X_t$ which is $K$-lipschitz.
\end{enumerate} 
\label{ItemLocalProduct}
\end{enumerate}
Each fiber $X_t$, $t \in T$, is called a \emph{horizontal leaf} in $X$. If
$L$ is a bi-infinite line in $T$ then $X_L$ is called a \emph{hyperplane}
in $X$. Items \ref{ItemNonincreasing} and \ref{ItemTreePathIsometry}
combine to show that the map of item \ref{ItemTreePathIsometry} is an
isometric embedding; the image $\Theta^\inv(x \cross T)$ is called a
\emph{vertical leaf} in $X$. For each subtree $T' \subset T$, the closest
point retraction $r \from T \to T'$ induces a map
$$X \xrightarrow{\Theta} F \cross T \xrightarrow{\Id \cross r} F \cross T'
\xrightarrow{\Theta^\inv} X_{T'}
$$ 
called \emph{vertical projection} of $X$ to $X_{T'}$.

\begin{Remark} Suppose $\Gamma$ is a graph of groups taken from a category
$\mathcal C$ as above. Let $M_\Gamma$ and $X_\Gamma \to T_\Gamma$ be as
constructed above starting from $\Gamma$. Then elementary constructions
produce a metric and a simplicial structure on $M_\Gamma$ which lifts to a
$\Gamma$-equivariant metric and simplicial structure on $X$ such that $X
\to T_\Gamma$ is a metric fibration. Item (\ref{ItemUnifContr}) follows
by compactness of $M_\Gamma$.
\end{Remark}

\begin{Remark} The definition has some redundancy: item
(\ref{ItemUnifContr}) is a formal consequence of item
(\ref{ItemHorizUnifContr}), as can be seen by elementary but mildly
tedious arguments. But by the previous remark we may dispense with these
arguments for the examples at hand.
\end{Remark}

The following lemma, applied to a bi-infinite line in $T$, gives good
geometric properties for hyperplanes:

\begin{lemma} If $\pi \from X \to T$ is a metric fibration then
there exist functions $\delta' \from [0,\infinity) \to [0,\infinity)$ and
$\rho
\from [0,\infinity)\to[0,\infinity)$
with $\displaystyle\lim_{t\to\infinity} \rho(t)=\infinity$, such that for
any subtree $T' \subset T$ we have:
\begin{itemize}
\item[(1)] The embedding $X_{T'} \to X$ is $\rho$-uniformly proper.
\item[(2)] The geodesic metric space $X_{T'}$ is $\delta'$-uniformly
contractible.
\end{itemize}
\label{lemma:metricfibration}
\end{lemma}

\begin{proof} To prove (1), consider $x,y \in X_{T'}$, let $D=d_X(x,y)$,
and let $\gamma\from[0,D]\to X$ be a geodesic connecting $x$ and $y$. Let
$N_D(T')$ be the $D$-neighborhood of $T'$ in $T$, so $\gamma \subset
X_{N_D(T')}$. Applying item \ref{ItemLocalProduct} iteratively,
projecting inward starting from the edges of $N_D(T')$ furthest
from $T'$, it follows that vertical projection $X_{N_D(T')}\to X_{T'}$
distorts any distance $r$ by at worst $K^D r$, and so
$d_{X_{T'}}(x,y) \le K^D D$.

To prove (2), suppose that $A \subset X_{T'}$ and $\diam_{X_{T'}}(A) \le
R$ and so $A$ is $R'$-homotopic to a constant in $X$ where $R'$ depends on
$R$ but not on $A$. This homotopy may then be mapped back to $X_{T'}$ by
vertical projection, distorting diameters of homotopy tracks by an amount
bounded in terms of $R'$ as we saw above. The result is an $R''$-homotopy
of $A$ to a constant in $T'$, with $R''$ depending only on $R$ and not on
$A$.
\end{proof}

Here is our generalization of Step 1. It applies to any metric fibration
of the form $X_\Gamma \to T_\Gamma$, where $\Gamma$ is a finite,
geometrically homogeneous graph of fundamental groups of manifolds in any
of the categories $\mathcal C$ described earlier.

\begin{theorem}
Let $\pi \from X \to T$ be a metric fibration whose fibers are
contractible $n$-manifolds for some $n$. Let $P$ be a contractible
$(n+1)$-manifold which is a uniformly contractible, bounded geometry,
metric simplicial complex. Then for any uniformly proper embedding $\phi
\from P \to X$, there exists a unique hyperplane $Q\subset X$ such that
$\phi(P)$ and $Q$ have finite Hausdorff distance in $X$. The bound on
Hausdorff distance depends only on the metric fibration data for $\pi$,
the uniform contractibility data and bounded geometry data for $P$, and
the uniform properness data for~$\phi$.
\label{theorem:metricfibration}
\end{theorem}

\begin{proof} Uniqueness of $Q$ follows obviously from the fact that
distinct hyperplanes in $X$ have infinite Hausdorff distance.

For existence of $Q$ we follow closely the proof of Proposition 4.1 of
\cite{FarbMosher:BSOne}, concentrating on details needed to explicate
the difference between the ``quasi-isometric'' setting of
\cite{FarbMosher:BSOne} and the present ``uniformly proper'' setting.

Using the bounded geometry of $P$, uniform contractibility of $X$, and
uniform properness of $\phi$, we may replace $\phi$ by a continuous,
uniformly proper map, moving values of $\phi$ a bounded distance.
Henceforth we shall assume $\phi$ is continuous.

Pick a topologically proper embedding of $T$ in an open disc $D$. For each
component $U$ of $D-T$, the frontier of $U$ in $D$ is
a bi-infinite line $L(U)$ in $T$. There is a homeomorphism of
pairs $(\overline U,L(U)) \approx (L(U)\cross [0,\infinity),L(U) \cross
0)$.

Consider the topologically proper embedding $X \xrightarrow{\Theta} F
\cross T\inject F\cross D$. Note that $F \cross D$ is a contractible
$(n+2)$-manifold. For each component $U$ of $D-T$ we have a
homeomorphism 
\begin{align*}
F \cross \overline U &\xrightarrow{\approx} F \cross (L(U) \cross
[0,\infinity))
 \xrightarrow{\approx} (F \cross L(U)) \cross [0,\infinity) \\
 & \xrightarrow[\Theta \cross \Id]{\approx} X_{L(U)} \cross [0,\infinity)
\end{align*}
The frontier of this set in $F \cross D$ is $F \cross L(U) \approx
X_{L(U)}$. Put a product metric and a product simplicial structure on
$X_{L(U)} \cross [0,\infinity)$ and glue to $F\cross L(U)$. Doing this
for each $U$, we impose a proper geodesic metric on $F \cross D$ for
which the inclusion $X\inject F \cross D$ is an isometric embedding.

The simplicial structure on $F \cross D$ evidently has bounded geometry.
Also, the metric space $F \cross D$ is uniformly contractible. To see
this, let $A\subset F\cross D$ have diameter $\le r$. If $A\intersect
X\ne\emptyset$ then homotoping along product lines of $X_{L(U)} \cross
[0,\infinity)$ for each $U$ we obtain an $r$-homotopy of $A$ into $F
\cross T \approx X$, and then we use uniform contractibility of $X$.
Whereas if $A\intersect X =\emptyset$, then $A \subset F \cross U \approx
X_{L(U)}\cross (0,\infinity)$ for some component $U$ of $D-T$; there is an
$r$-homotopy of $A$ into some $X_{L(U)}\cross x$, and the latter is
uniformly contractible by Lemma \ref{lemma:metricfibration}.

We now plug this setup into the coarse separation and packing methods of
Farb--Schwartz \cite{FarbSchwartz} and Schwartz
\cite{Schwartz:Diophantine}. We'll use a generalization of the Coarse
Separation Theorem with more easily applied hypotheses, due to
Kapovich--Kleiner \cite{KapovichKleiner:duality}. We denote the
\nb{$r$}ball about a subset $A$ of a metric space $M$ by $B_r(A;M)$. In a
metric space $Z$, a subset $U\subset Z$ is \emph{deep in $Z$} if for each
$r>0$ there exists $x \in U$ such that $B_r(x;Z)\subset U$. A subset $A
\subset Z$ \emph{coarsely separates} $Z$ if for some $D>0$ there are at
least two components of $Z-N_D(A;Z)$ which are deep in $Z$; the constant
$D$ is called a \emph{coarse separation constant} for $A$. Note that if
subsets $A$ and $B$ of $Z$ have bounded Hausdorff distance from each
other, then $A$ coarsely separates $Z$ if and only if $B$ does. 

Here is an elementary consequence of the definitions:

\begin{lemma}
\label{lemma:sep2}
Let $f\from X\to Y$ be a quasi-isometry between geodesic metric spaces.
If $A \subset X$ coarsely separates $X$ then $f(A)$ coarsely separates
$Y$, with separation constant depending only on the quasi-isometry 
constants of $f$ and the separation constant for $A$.
\qed\end{lemma}

Here is the version of the Coarse Separation Theorem that we will use.

\begin{theorem}[\cite{KapovichKleiner:duality}]
Let $P$ be a contractible $(n+1)$-manifold, $Z$ a contractible
$(n+2)$-manifold, and suppose that $P,Z$ are uniformly contractible,
bounded geometry, metric simplicial complexes. Let $\Phi \from P \to Z$
be a uniformly proper map. Then $\Phi(P)$ coarsely separates
$Z$, with coarse separation constant $D$ depending only on the uniform
contractibility and bounded geometry data for $P$ and $Z$ and the uniform
properness data for $\Phi$. Moreover if $\Phi$ is continuous then we may
take $D=0$, that is, $Z-\Phi(P)$ has at least two components which are
deep in $Z$.
\label{theorem:coarse:separation}
\qed\end{theorem}

\begin{Remark} In fact there are exactly two components of
$Z-N_D(\Phi(P);Z)$ which are deep in $Z$ (see
\cite{KapovichKleiner:duality}).
\end{Remark}

Following \cite{FarbSchwartz} we have a corollary:

\begin{theorem}[Packing Theorem] Let $Q,P$ be contractible
$(n+1)$-manifolds, which are uniformly contractible, bounded geometry,
metric simplicial complexes. Let $\psi \from Q \to P$ be a
uniformly proper map. Then there exists $R>0$ such that $N_R(\psi(Q);P) =
P$. The constant $R$ depends only on the uniform contractibility data and
bounded geometry data for $Q,P$ and the uniform properness data for
$\psi$.
\label{theorem:packing}
\end{theorem}

\begin{proof}
If no such $R$ exists then the image of the map $Q \xrightarrow{\psi} P
\inject P\cross\R$ does not coarsely separate $P \cross \R$, violating 
Theorem \ref{theorem:coarse:separation}.
\end{proof}

Continuing with the proof of Theorem \ref{theorem:metricfibration},
compose the continuous, uniformly proper map $\phi \from P \to X$ with
the isometric embedding $X \to F \cross D$ to obtain a
continuous, uniformly proper map $\Phi \from P\to F \cross D$. By
the Coarse Separation Theorem it follows that $(F \cross D) - \Phi(P)$
has at least two components which are deep in $F \cross D$.

Now take the argument of \cite{FarbMosher:BSOne}, Step 1, pages 426--427,
and apply it verbatim, to produce a hyperplane $Q \subset X$ such that $Q
\subset \Phi(P)$. Next take the argument of Step 2, pages 427--428, and
apply it verbatim, replacing ``quasi-isometric embeddings'' with
``uniformly proper maps'' and using the Packing Theorem above, to show
the existence of $R'$ such that $\phi(P) \subset N_{R'}(Q;X)$, where $R'$
depends only on the metric fibration data for $\pi$, the uniform
contractibility and bounded geometry data for $P$, and the uniform
properness data for $\phi$. 

This finishes the proof of Theorem \ref{theorem:metricfibration} and of
Step~1.
\end{proof}

\paragraph*{Step 2. A quasi-isometry takes hyperplanes and horizontal
leaves in $X_M$ to hyperplanes and horizontal leaves in $X_N$:}
\hfill\break\noindent
Consider integral matrices $M,N \in \GL_\cross(n,\R)$ with $\det M,
\det N>1$, and let $f\from X_M \to X_N$ be a quasi-isometric embedding.
Then there is a constant $A \ge 0$, depending only on $X_M, X_N$ and the
quasi-isometry constants of $f$, such that:
\begin{enumerate}
\item For each hyperplane $P \subset X_M$ there exists a unique
hyperplane $Q \subset X_N$ such that
$d_\Haus(f(P),Q) \le A$.
\item For each horizontal leaf $L$ of $X_M$ there exists a horizontal
leaf $L'$ of $X_N$ such that
$d_\Haus(f(L),L')\leq A$.
\end{enumerate}

The proof of this step is the first place in our arguments where the
assumption that $\det M,\det N>1$ is crucial. Again we will investigate
this step in the general setting of metric fibrations over trees. 

Consider a metric fibration $\pi \from X \to T$. The tree $T$ is
\emph{bushy} if there exists a constant $\beta$ such that each point of
$T$ is within distance $\beta$ of some vertex $v$ such that
$T-v$ has at least 3 unbounded components. Note that if $M$ is an integer
matrix in $\GL_\cross(n,\R)$,  and if $X_M \to T_M$ is the associated
metric fibration over the Bass-Serre tree $T_M$ of the group $\Gamma_M$,
then $T_M$ is bushy if and only if $\det M > 1$. In fact, for any graph
of finitely generated groups, the Bass-Serre tree is either bounded,
quasi-isometric to a line, or bushy, and the question of which
alternative holds is easily decided by inspection of the graph of groups.

Here is our generalization of Step 2:

\begin{theorem}
Let $\pi \from X \to T$, $\pi' \from X' \to T'$ be  metric fibrations
over $\beta$-bushy trees $T,T'$, such that the fibers of $\pi$ and $\pi'$
are contractible $n$-manifolds for some $n$. Let $f\from X \to X'$ be a
quasi-isometry. Then there exists a constant $A$, depending only on the
metric fibration data of $\pi,\pi'$, the quasi-isometry data for $f$, and
the constant $\beta$, such that:
\begin{enumerate}
\item[(1)] For each hyperplane $P \subset X$ there exists a unique
hyperplane $Q \subset X'$ such that $d_\Haus(f(P),Q) \le A$.
\item[(2)] For each horizontal leaf $L \subset
X$ there is a horizontal leaf $L' \subset X'$ such that $d_\Haus(f(L),L')
\le A$.
\end{enumerate}
\label{theorem:horizontal}
\end{theorem}

\begin{Remark}
This result is used in \cite{MosherSageevWhyte} to prove quasi-isometric
ridigity for fundamental groups of geometrically homogeneous graphs of
groups whose vertex groups are fundamental groups of manifolds in a
category $\mathcal C$ as above, as long as that class of groups is itself
quasi-isometrically rigid. For example, quasi-isometric rigidity is
proved for graphs of $\Z$'s, $\Z^n$'s, surface groups, lattices in
semisimple Lie groups, nilpotent groups, etc.
\end{Remark}

\begin{proof} To prove (1), by Lemma \ref{lemma:metricfibration} the
inclusion map $P\inject X$ is uniformly proper and $P$ is uniformly
contractible, and clearly $P$ is a contractible $(n+1)$-manifold.
Composing with $f$ we obtain a uniformly proper map $P \to X'$. Now apply
Theorem
\ref{theorem:metricfibration}.

The idea of the proof of (2) is that bushiness of the tree allows one to
gain quasi-isometric control over horizontal leaves by considering them
as ``coarse intersections'' of hyperplanes.

\begin{Definition}[Coarse intersection]
A subset $W$ of a metric space $X$ is a {\em coarse intersection}
of subsets $U,V \subset X$, denoted $W=U\coarsecap V$, if there
exists $K_0$ such that for every $K\geq K_0$ there exists $K'\geq
0$ so that 
$$d_\Haus\bigl(\nbhd_K(U)\cap \nbhd_K(V),W\bigr)\leq K'
$$  
Note that although such a set $W$ may not exist, when it does
exist then any two such sets are a bounded Hausdorff distance from
each other.
\end{Definition}

The following fact is an elementary consequence of
the definitions.

\begin{lemma}
\label{lemma:intersection2}
For any quasi-isometry $f\from X\to Y$ of metric spaces, and $U,V \subset 
X$, if $U \coarsecap V$ exists then $f(U \coarsecap V)$ is a coarse 
intersection of $f(U)$, $f(V)$, with constants depending only on the 
quasi-isometry constants for $f$ and the coarse intersection constants for 
$U$ and $V$.
\qed\end{lemma}

Consider now a metric fibration $\pi \from X \to T$. A subset of
$X$ of the form $X_\sigma = \pi^\inv(\sigma)$, where $\sigma$
is an infinite ray in $T$, will be called a {\em half-plane} in
$X$. The next lemma is an easy observation---see
Figure~\ref{FigureCoarseInt}.

\begin{figure}
\centeredepsfbox{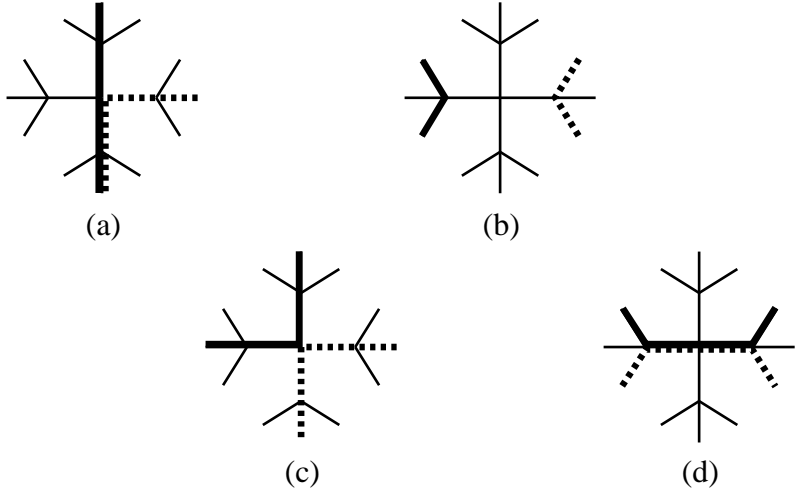}
\caption{Possible coarse intersections of distinct hyperplanes in
$X$, projected to $T$. In (a), $P_1 \coarsecap P_2 = P_1
\intersect P_2$ is a half-plane. In (b--d), $P_1 \coarsecap P_2$
is a horizontal leaf; $P_1 \intersect P_2$ can be empty (b), a
horizontal leaf (c), or a finite strip of horizontal leaves (d).}
\label{FigureCoarseInt}
\end{figure}

\begin{lemma}
\label{lemma:intersection1}
Let $\pi \from X \to T$ be a metric fibration over a tree $T$.
Let $P_1$ and $P_2$ be distinct hyperplanes in $X$.  Then  $P_1\coarsecap
P_2$ exists and is a bounded Hausdorff
distance from either a half-plane or a horizontal leaf in $X$. Moreover,
$P_1\coarsecap P_2$ is a bounded Hausdorff distance from a half-plane if
and only if $P_1\cap P_2$ is a half-plane.
\qed\end{lemma}

We remark that $P_1\coarsecap P_2$ can be an arbitrarily large
finite Hausdorff distance from a horizontal leaf; see Figure 
\ref{FigureCoarseInt}(b,d).

\begin{lemma}
\label{lemma:intersection3}
Let $\pi \from X \to T$, $\pi' \from X' \to T'$ be metric
fibrations. Let $f \from X \to X'$ be a quasi-isometry. Suppose $P_1$ and
$P_2$ are distinct hyperplanes in $X$ which intersect in a half-plane.
Then $f(P_1)$ and $f(P_2)$ are a uniformly bounded Hausdorff distance from
distinct hyperplanes $Q_1,Q_2$ in $X'$ which intersect in a half-plane in
$X'$.
\end{lemma}

\begin{proof}
By Theorem \ref{theorem:metricfibration}, there exists a constant $A$ so
that $f(P_i)$ is within Hausdorff distance $A$ of a unique hyperplane
$Q_i$ in $X'$. Since $P_1, P_2$ are distinct they have infinite Hausdorff
distance, so $Q_1$ and $Q_2$ have infinite Hausdorff distance and hence
$Q_1 \ne Q_2$. 

By Lemma \ref{lemma:intersection1}, it is enough to prove
that $Q_1\coarsecap Q_2$ is not a bounded Hausdorff distance from a
horizontal leaf in $X'$. If $Q_1\coarsecap Q_2$ is a bounded Hausdorff
distance from a horizontal leaf, then since any horizontal leaf in $Q_1$
coarsely separates $Q_1$ it must be that $Q_1\coarsecap Q_2$ coarsely
separates $Q_1$. But $P_1\coarsecap P_2$ does not coarsely
separate $P_1$. This contradicts Lemma \ref{lemma:sep2}.
\end{proof}

We now prove Theorem \ref{theorem:horizontal}. Consider the
quasi-isometry $f \from X \to X'$. Since $T$ is bushy, any horizontal
leaf $L$ in $X$ can be realized as a coarse intersection of three
hyperplanes $P_1, P_2, P_3$, such that the pairwise intersections $P_1
\intersect P_2$, $P_2 \intersect P_3$, $P_3 \intersect P_1$ form three
half-planes, any two of which have infinite Hausdorff distance.
Moreover, $d_\Haus(L,P_1\cap P_2 \cap P_3) \le \beta$ where
$\beta$ is a bushiness constant for $T$ (see Figure \ref{FigureTriple}).

Consider the unique hyperplane $Q_i$ which lies a Hausdorff
distance of at most $A$ from $f(P_i)$, $i=1,2,3$. By Lemma
\ref{lemma:intersection3}, the pairwise intersections $Q_1 \intersect
Q_2$, $Q_2 \intersect Q_3$, $Q_3 \intersect Q_1$ are all half-planes, any
two of which have infinite Hausdorff distance. The following elementary
fact about trees, applied to $T'$, now shows that $Q_1\cap Q_2\cap Q_3$ is
a horizontal leaf $L'$ in $X'$:

\begin{figure}
\centeredepsfbox{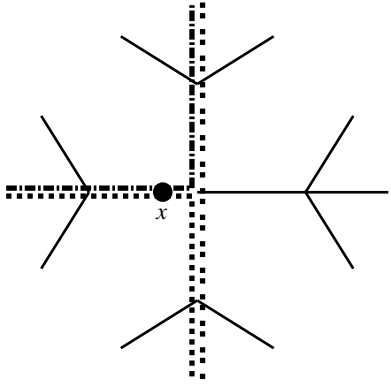}
\caption{Any point $x \in T$ is a bounded distance $\beta$ from a vertex
$v \in T$ that separates $T$ into at least three unbounded components.
The vertex $v$ is the (coarse) intersection of three proper lines
$\ell_1,\ell_2,\ell_3$, such that the pairwise intersections $\ell_1
\intersect\ell_2$, $\ell_2\intersect \ell_3$, $\ell_3 \intersect \ell_1$
are rays in $T$, any two of which have infinite Hausdorff distance.
Moreover, $d(x,\ell_1 \cap\ell_2\cap \ell_3) \le\beta$.}
\label{FigureTriple}
\end{figure} 

\medskip \noindent {\bf Fact about trees: }
Let $\ell_1,\ell_2,\ell_3$ be bi-infinite lines in a simplicial tree $T'$,
such that the pairwise intersections $\ell_1 \intersect \ell_2$, $\ell_2
\intersect
\ell_3$, $\ell_3 \intersect \ell_1$ are all infinite rays in $T'$, any two
of which have infinite Hausdorff distance. Then $\ell_1\cap \ell_2\cap
\ell_3$ is a vertex of $T'$.
\medskip

Since $L \subset N_\beta(P_i)$ it follows that 
$$f(L) \subset N_{K\beta+C}(f(P_i)) \subset N_{K\beta+C+A}(Q_i), \quad
i=1,2,3
$$ 
But clearly we have $\bigcap_{i=1}^3 N_{K\beta+C+A}(Q_i) =
N_{K\beta+C+A}(L')$. 

To summarize, given a horizontal leaf $L$ of $X$, we have found a 
horizontal leaf $L'$ of $X'$ such that $L \subset N_{A'}(L')$ where
$A'=K\beta+C+A$. A similar argument using a coarse inverse for $f$
provides the desired bound for $d_\Haus(f(L),L')$. This completes the
proofs of Theorem \ref{theorem:horizontal} and of Step 2. \end{proof}

\paragraph*{Step 3. A quasi-isometry takes coherent hyperplanes in
$X_M$ to coherent hyperplanes in $X_N$.}
\hfill\medskip 

Let $M,N$ be as in the statement of Proposition
\ref{proposition:induced}, and fix a quasi-isometry $f \from X_M \to
X_N$.

Let $P$ be any coherent hyperplane in $X_M$.  By Step 2 it follows that
$f(P)$ is within a Hausdorff distance $A$ from a unique hyperplane $Q$ in
$X_N$.  By composing $f \restrict P$ with vertical projection $X_N \to Q$
we obtain a map $\phi \from P\to Q$. The inclusion maps $P \inject X_M$
and $Q \inject X_N$ are coarsely lipschitz and uniformly proper; indeed
they are isometric embeddings with respect to the induced path metrics on
$P,Q$. By Lemma \ref{lemma:path:subspace}, $\phi$ is a quasi-isometry,
with quasi-isometry constants depending only on those for $f$. By Step 2,
$f$ coarsely respects the horizontal foliations of $X_M$ and $X_N$;
vertical projection $X_N \to Q$ takes horizontal leaves to horizontal
leaves, and so $\phi$ coarsely respects the horizontal foliations of $P$
and $Q$, with a coarseness constant depending only on the quasi-isometry
constants of $f$.
  
Since $P$ is a coherent hyperplane it is isometric to $G_M$. Since
$Q$ is a hyperplane it is isometric to either $G_N$ or $H_N$, and
we now show that the second possibility cannot occur.

\begin{proposition}
Given matrices $M,N \in \GL_\cross(n,\R)$ with 
$\det M, \det N > 1$, there is no quasi-isometry $\phi\from G_M\to H_N$
which coarsely respects horizontal foliations.
\label{proposition:area} \end{proposition}

\begin{proof} The idea of the proof is to compare the growth types of the 
filling area functions for ``quasivertical bigons'' in $G_M$ and in $H_N$. 
In $G_M$ this growth type will be quadratic, while in $H_N$ it will be 
exponential.

Let $H = G_M$, $H_M$, $G_N$, or $H_N$. There is a quotient map $H \to \R$
whose point pre-images give  the horizontal foliation of $H$, and such
that the Hausdorff distance  between two horizontal leaves equals the
distance between the  corresponding points in $\R$. A path $\gamma$ in
$H$ is said to be  $(K,C)$-\emph{quasivertical} if its projection to $\R$
is a $(K,C)$-quasigeodesic. Define a $(K,C)$-\emph{quasivertical bigon}
in $H$ to  be a pair of $(K,C)$-quasivertical paths $\gamma,\gamma'$
which begin and  end at the same point. 

If $K,C$ are fixed, we define a filling area function $A(L)$ for
$(K,C)$-quasivertical bigons in $H$. Given a $(K,C)$-quasivertical
bigon $\gamma,\gamma'$, its \emph{filling area} is the infimal
area of a Lipschitz map $D^2 \to H$ whose boundary is a
reparameterization of the closed curve $\gamma^\inverse *
\gamma'$; such a map $D^2 \to H$ is called a \emph{filling
disc} for $\gamma^\inverse * \gamma'$. For each $L \ge 0$ define
$\A(L)$ to be the supremal filling area over all $(K,C)$-quasivertical
bigons $\gamma,\gamma'$ in $H$ such that $\Length(\gamma) +
\Length(\gamma') \le L$.

Suppose there is a quasi-isometry $\phi \from G_M \to H_N$ which coarsely
respects horizontal foliations. Let $\bar\phi \from H_N \to G_M$ be a
coarse inverse for $\phi$, also coarsely respecting horizontal foliations.
Clearly $\bar\phi$ takes any $K,C$-quasivertical bigon in
$H_N$ to a $(K',C')$-quasivertical bigon in $G_M$, distorting lengths by
at worst an affine function; this affine function, and the constants
$K',C'$, depend only on $K,C$, the quasi-isometry constants for $\phi$,
and the Hausdorff constant for the induced height function. Fill the
resulting bigon in $G_M$ as efficiently as possible, and map back to
$H_N$ via $\phi$, distorting area by at worst an affine function which
again has the same dependencies. We thereby obtain a filling of the
original bigon in $H_N$. If $\A_1(L)$ denotes the filling area function
for
$(K',C')$-quasivertical bigons in $G_M$, and if $\A_2(L)$ denotes the
filling area function for $(K,C)$-quasivertical bigons in $H_N$, it
follows that the growth type of $\A_2(L)$ is dominated by the growth type
of $\A_1(L)$, that is,
$$\A_2(L) \le \alpha \cdot \A_1(\beta L + \delta) + \zeta
$$
for some positive constants $\alpha,\beta,\delta,\zeta$ independent of
$L$.

However, we shall now show that $\A_1(L)$ has a quadratic upper bound while 
$\A_2(L)$ has an exponential lower bound, contradicting the above 
inequality.

Consider a $K',C'$-quasivertical bigon $\gamma,\gamma'$ in $G_M$.
Applying the argument of Claim \ref{claim:quasivertical}, there are
center leaves $\tau,\tau'$ in $G_M$ and quasivertical paths $\rho
\subset \tau,\rho' \subset \tau'$ which stay uniformly close to
$\gamma,\gamma'$, respectively. The initial points of $\rho,\rho'$
are at a uniformly bounded distance, as are the terminal points,
and it follows that $\rho'$ stays uniformly close to a
quasivertical path $\rho'' \subset \tau$. Connecting initial and
terminal endpoints with short paths $\eta,\eta'$ we thus obtain a
closed curve $\rho^\inverse * \eta * \rho'' * \eta'$, contained in
a center leaf of $G_M$, which stays uniformly close to
$\gamma^\inverse * \gamma'$. Since center leaves of $G_M$ are
isometric to Euclidean space, in which the filling function is
quadratic, it follows that $\A_1(L)$ has a quadratic upper bound.

To show that $\A_2(L)$ has an exponential lower bound, we now
construct quasivertical bigons in $H_N$ which can be filled only
by discs of exponential area. In the case where $N$ is a $1 \cross
1$ matrix such loops are given explicitly in \cite{ECHLPT}, Chapter 7.4; 
examples for general $N$ are
simple modifications of this example. To be explicit, choose an
eigenvalue of $N$ of absolute value $\alpha>1$; such an eigenvalue
exists because $\det N > 1$. Choose an affine subspace $A \subset
\R^n$ parallel to the $\alpha$-eigenspace of $N$. Consider the
subspace $A \cross \reals \subset \reals^n \cross \reals \approx
G_N$.

For each fixed $L \ge 0$, choose two vertical segments 
$g,g'$ in $A \cross [0,\infinity)$ whose upper endpoints are in $A \cross 
L$ and whose lower endpoints are in $A \cross 0$, and so that the distance 
in $A \cross L$ between the upper endpoints, measured using the Riemannian 
metric on $G_N$, is equal to $1$; it follows that the distance in $A 
\cross 0$ between the lower endpoints, measured using the Riemannian 
metric on $G_N$, is within a constant multiple of $\alpha^L$. 

Now double 
this picture, in the doubled $G_N$ horoball $H_N$, to get a closed loop in 
$H_N$, that is: in one horoball go up $g$, across $1$ unit, and down $g'$, 
and then in the other horoball go up $g'$, across $1$ unit, and down $g$; 
let $\rho$ be the resulting closed curve in $H_N$. We have $\Length(\rho) 
= 4L+2$. To see that the filling area of $\rho$ is exponential in $L$, 
note that any filling disc for $\rho$ must contain a path in $A \cross 0$ 
connecting the lower endpoints of $g,g'$, because $A \cross 0$ separates 
the two halves of $\rho$ in $H_N$. This path has length exponential in 
$L$; and a neighborhood of this path in the filling disc has area 
exponential in $L$.
\end{proof}

\subsubsection*{Step 4. A horizontal respecting quasi-isometry preserves
transverse orientation}
Let $M$, $N$, and $f \from \Gamma_M \to \Gamma_N$ be as in the statement
of Proposition \ref{proposition:induced}. By Step 3 there is a
quasi-isometry $\phi\from G_M\to G_N$, and by Step 2, $\phi$ coarsely
respects the horizontal foliations of $G_M$ and $G_N$. Suppose that
$\phi$ reverses the transverse orientation. There  is a quasi-isometry
$G_N \to G_{N^\inverse}$ which coarsely respects  horizontal foliations,
\emph{reversing} transverse orientations.  Precomposing with $\phi \from
G_M \to G_N$ and applying Steps 1--3, we  obtain a quasi-isometry $G_M
\to G_{N^\inverse}$ which coarsely respects  the transversely oriented
horizontal foliations. Applying Theorem
\ref{theorem:horizontal:respecting}, it follows that $M$ and
$N^\inverse$  have positive real powers with the same absolute Jordan
form, and so these  powers also have the same determinant. But each 
positive power of $M$ has determinant $>1$, whereas every  positive power
of $N^\inverse$ has determinant $<1$, a  contradiction showing that $\phi$
must preserve the transverse orientation.

This completes the proof of Proposition \ref{proposition:induced}.\qed

\begin{Remark}
Note in the proof of Proposition \ref{proposition:induced} 
that different choices of coherent hyperplanes
in $X_M$ yield different quasi-isometries $\phi$. In some cases
$\phi$ is well-defined up to some constant $A$, that is, for any
two choices of coherent hyperplane in $X_M$, the induced maps
$\phi_1,\phi_2 \from G_M \to G_N$ satisfy $\sup_x
d(\phi_1(x),\phi_2(x)) \le A$. This is true, for example, in the
``centerless'' case where $M,N$ have no eigenvalues on the unit
circle. In the general case, the best that can be said is that the
map induced by $\phi$ from the center leaf space of $G_M$ to the
center leaf space of $G_N$ is well defined up to a constant, with
respect to the Hausdorff metrics on the center leaf spaces.
\end{Remark}


\section{Finding the Integers}
\label{section:classification}

In this section we prove Theorem \ref{theorem:classification}. 
Let $M,N$ be integral $(n \cross n)$ matrices with 
$\abs{\det M}, \abs{\det N} > 1$. We must prove that $\Gamma_M$ 
is quasi-isometric to $\Gamma_N$ if and only if there exist 
positive integers $a,b$ such that $M^a$ and $N^b$ have the 
same absolute Jordan form.

First we show that the groups $\Gamma_{M^a}$ and $\Gamma_M$ are
quasi-isometric, for any positive integer $a$, by showing that
$\Gamma_{M^a}$ is a subgroup of finite index in $\Gamma_M$,
specifically of index $a$. To see why, consider the presentations
\begin{align*}
\Gamma_M &= \left< \Z^n,t \suchthat t^\inverse x t = M(x), x \in \Z^n
\right> \\
\Gamma_{M^a} &= \left< \Z^n,s \suchthat s^\inverse x s = M^a(x), x \in \Z^n \right>
\end{align*}
Define a homomorphism $\Gamma_M \mapsto \Z / a\Z$ by $\Z^n \mapsto
0, t \mapsto 1$. This homomorphism is onto, and its kernel is
generated by $\Z^n,t^a$. This kernel is isomorphic to
$\Gamma_{M^a}$ under the injection $\Gamma_{M^a} \inject \Gamma_M$
given by $x \mapsto x, s \mapsto t^a$.

Similarly, $\Gamma_{N^b}$ is quasi-isometric to $\Gamma_N$, for any
positive integer $b$.

By squaring $M,N$ if necessary, we may therefore assume that 
$\det M, \det N > 1$, and that $M$ and $N$ lie on \nb{1}parameter 
subgroups; we continue with this assumption up 
through the end of the proof in \S\ref{section:secondhalf}. 
Choose \nb{1}parameter subgroups $M^t,N^t$ of $\GL(n,\R)$
with $M=M^1, N=N^1$, let $G_M,G_N$ be the associated Lie 
groups constructed in \S\ref{section:liegroup}, and let 
$X_M,X_N$ be the associated geodesic metric spaces 
constructed in \S\ref{section:coarsetop}. The group 
$\Gamma_M$ is quasi-isometric to $X_M$, and 
$\Gamma_N$ is quasi-isometric to $X_N$. 

\subsection{The first half of the classification}

Assuming that $M^a$ and $N^b$ have the same absolute Jordan form,
where $a,b$ are positive integers, we must prove that $\Gamma_M$
and $\Gamma_N$ are quasi-isometric. We have shown above that
$\Gamma_{M^a}$ and $\Gamma_M$ are quasi-isometric, and that
$\Gamma_{N^b}$ and $\Gamma_N$ are quasi-isometric. Replacing $M$
by $M^a$ and $N$ by $N^b$, we may therefore assume that $M,N$ have
the same absolute Jordan form. We shall prove that $\Gamma_M,
\Gamma_N$ are quasi-isometric by constructing a bilipschitz
homeomorphism between $X_M$ and $X_N$.

Since the absolute Jordan forms of $M,N$ are equal it follows that
$\det M = \det N$; let $d$ be the common value. Applying Proposition
\ref{proposition:choices}, there is a bilipschitz homeomorphism from $G_M
= \R^n \semidirect_M \R$ to $G_N = \R^n\semidirect_M \R$ of the form
$(x,t) \mapsto (Ax,t)$ for some $A\in \GL(n,\R)$. In the fiber product
description of $X_M$, $X_N$, the trees $T_M$ and $T_N$ may both be
identified with the homogeneous, oriented tree $T_d$ with one incoming
and $d$ outgoing edges at each vertex. The bilipschitz homeomorphism
$G_M\to G_N$ and the identity homeomorphism $T_d \to T_d$ both respect
the height functions, and so these two homeomorphisms combine to give the
desired bilipschitz homeomorphism $X_M \to X_N$.

\subsection{Quasi-isometric implies integral 
powers have same absolute Jordan forms} 
\label{section:secondhalf}

Assuming $\Gamma_M, \Gamma_N$ are quasi-isometric, there is a quasi-isometry $f \from X_M \to X_N$. Combining Proposition \ref{proposition:vhpreserved} and Theorem
\ref{theorem:horizontal:respecting} gives $r\in\R_+$ such that
$M^r$ and $N$ have the same absolute Jordan form.  We must show
that there exist $a,b\in \Z_+$ so that $M^a$ and $N^b$ have the
same absolute Jordan form.

Since $M^r$ and $N$ have the same absolute Jordan form, listing
the absolute values of the eigenvalues of $M$ and $N$ in
increasing order we obtain 
\begin{align*}
\mu_{-a} \le \cdots \le \mu_{0} \le 1 < \mu_{1}=:\alpha_M \le \cdots \le \mu_b \\
\nu_{-a} \le \cdots \le \nu_{0} \le 1 < \nu_{1}=:\alpha_N \le \cdots \le \nu_b
\end{align*}
with $\mu_i^r = \nu_i$, $-a \le i \le b$. From this it follows that
$$\frac{\log\alpha_N}{\log\alpha_M} = r = \frac{\log\det
N}{\log\det M}
$$

Let $\Q_M$ denote the set of coherent hyperplanes in $X_M$, and let $h_M$ 
denote the height function on $M$.  We define a metric on $\Q_M$ as 
follows: given coherent hyperplanes $P_1,P_2$, let $L$ denote 
the horizontal leaf $L=\partial(P_1\cap P_2)$. Then we set 
$$d_{\Q_M}(P_1,P_2)=\bigl(\det M\bigr)^{-h_M(L)}$$

It is easy to check that this defines a metric on $\Q_M$, and
since the tree $T_M$ branches $m=\det M$ times as $h_M$
increases by $1$, the metric space $(\Q_M,d_{\Q_M})$ is isometric
to the $m$-adic rational numbers in their usual metric of Hausdorff
dimension~1. Similarly, attached to $X_N$ is a metric space
$(\Q_N,d_{\Q_N})$ isometric to the $n$-adic rational numbers, with
$n=\det N$.

From Step 3 in the proof of Proposition
\ref{proposition:vhpreserved} (see \S\ref{section:coherent}), the
quasi-isometry $f:X_M\to X_N$ takes each coherent hyperplane in
$X_M$ to within a uniform Hausdorff distance of a unique coherent
hyperplane in $X_N$, hence induces a bijection $\psi\from \Q_M \to
\Q_N$. For each $\ell \in \Q_M$, setting $\ell' = \psi(\ell)$,
there is an induced horizontal-respecting quasi-isometry $P_\ell
\to P'_{\ell'}$, and by Time Rigidity (Proposition
\ref{proposition:height:rigidity}) this quasi-isometry has an
induced time change of the form $t \mapsto mt+b$ where
$$m = \frac{\log \alpha_M}{\log \alpha_N}=1/r\
$$
and where $b$ depends ostensibly on $\ell$. However, for another
$\ell_1$, $P_\ell$ and $P_{\ell_1}$ coincide below some value of $t$, and so $t \mapsto mt+b$ is an induced time change for both
$P_\ell \mapsto P'_{\ell'}$ and $P_{\ell_1} \mapsto P'_{\ell'_1}$,
possibly with a larger coarseness constant (this argument is taken
from Claim 6.3 on p.\ 436 of \cite{FarbMosher:BSOne}). Therefore,
there is a uniform induced time change $t \mapsto mt+b$ with $b$
independent of $\ell$, and with a uniform Hausdorff constant $A$.

We now claim that $\psi$ is a bilipschitz homeomorphism.  To this
end, let $P_1,P_2\in \Q_M$ be given. Let $L=\partial(P_1\cap
P_2)$ and let $L'=\partial(\psi(P_1)\cap \psi(P_2))$. Then
$$h_N(L') \ge m\cdot h_M(L) + b - A
$$
Hence 
\begin{align*}
\frac{ d_{\Q_N}(\psi(P_1),\psi(P_2))}
{ d_{\Q_M}(P_1,P_2)}
&=\frac{ (\det N)^{-h_N(L')}}
{ (\det M)^{-h_M(L)}}\\
&\leq \frac{ (\det N)^{-m h_M(L)-b+A}}
{ (\det M)^{-h_M(L)}}\\
&=\frac{\bigl((\det N)^{(\log\det M/\log\det N))}\bigr)^{-h_M(L)} (\det N)^{-b+A}}
{ (\det M)^{-h_M(L)}}\\
&=(\det N)^{-b+A}
\end{align*}
which is a constant not depending on $P_1$ or $P_2$.  Hence $\psi$
is Lipschitz.  The same argument applied to $\psi^{-1}$ shows that
$\psi$ is bilipschitz.

Applying Cooper's Theorem (appendix to \cite{FarbMosher:BSOne},
Corollary 10.11) on bilipschitz homeomorphisms of Cantor sets, we
obtain that there exist integers $a,b>0$ such that $(\det M)^a =
(\det N)^b$. Since $M^r$ and $N$ have the same absolute Jordan
form, we have
$$\frac{b}{a} = \frac{\log \det M}{\log \det N} = r
$$
and so $(M^r)^a = M^b$ and $N^a$ have the same absolute Jordan form.


\section{Quasi-isometric rigidity}
\label{section:qirigidity}

In this section we prove Theorem \ref{theorem:rigidity} in a series of
steps. Recall the hypotheses: $M$ is an integer matrix in $\GL(n,\R)$
with $\abs{\det M} > 1$, and $G$ is a finitely generated group
quasi-isometric to $\Gamma_M$. By squaring $M$ if necessary we may assume
that $M \in \GL_\cross(n,\R)$ and $\det M > 1$, and therefore $\Gamma_M$
is quasi-isometric to $X_M$. It follows that $G$ is quasi-isometric to
$X_M$.

\paragraph*{Step 1.} 
The action of $G$ on itself by left multiplication
can be conjugated by the quasi-isometry $G \to X_M$ to give a proper,
cobounded quasi-action of $G$ on $X_M$ (see \cite{FarbMosher:BSTwo},
Proposition 2.1). Since $\det M > 1$ we may apply Theorem
\ref{theorem:horizontal}, concluding that the quasi-action of $G$
on $X_M$ coarsely respects the fibers of the uniform metric fibration
$X_M \to T_M$.

\paragraph*{Step 2.} Now we use the following result of
\cite{MosherSageevWhyte}. Suppose $\pi \from X \to T$ is a uniform metric
fibration over a bushy tree $T$. If $G$ is a finitely presented group
with a cobounded, proper quasi-action on $X$, and if the quasi-action
coarsely respects the fibers, then $G$ is the fundamental group of a
graph of groups whose vertex and edge groups are quasi-isometric to a
fiber $X_t = \pi^\inv(t)$. 

By Step 1, this result applies to the quasi-action of $G$ on $X_M$,
because $G$ is quasi-isometric to the finitely presented group $\Gamma_M$
and so $G$ is finitely presented. The fibers of the map $X_M \to T_M$ are
isometric to $\R^n$, and it follows that $G$ is the fundamental group of a
graph of groups with each vertex and edge group quasi-isometric to $\R^n$.

\paragraph*{Step 3.} Any finitely-generated group quasi-isometric to
$\R^n$ is virtually $\Z^n$ (see \cite{Gersten:isoperimetric}), and so $G$
is the fundamental group of a graph of groups whose vertex and edge
groups are virtually $\Z^n$.

\paragraph*{Step 4.} Applying the argument in Section 5 of 
\cite{FarbMosher:BSTwo} to $G$ gives that either $G$ contains a noncyclic
free group or $G$ is an ascending HNN extension of the form
$$G = A_\phi = \left< A,t \suchthat tat^\inv = \phi(a), \forall a \in
A\right>
$$
where $A$ is virtually $\Z^n$ and $\phi \from A \to A$ is an injective
endomorphism.  Since $\Gamma_M$ is amenable, and since $G$ is
quasi-isometric to $\Gamma_M$, then $G$ is amenable, and so $G$
cannot contain a noncyclic free group. The second possibility
must therefore occur: $G = A_\phi$ as above.

\paragraph*{Step 5.}
Now we turn to an analysis of injective endomorphisms of virtually
abelian groups. Suppose $A$ is a finitely generated, virtually abelian
group. Any injective endomorphism of $A$ has finite index image. 

A subgroup $B\subset A$ is \emph{characteristic for endomorphisms} if,
for any injective endomorphism $\phi \from A \to A$, we have $\phi(B)
\subset B$.

Given a group $A$ and $g \in A$, the centralizer of $g$ in $A$ is denoted
$C_A(g)$. The \emph{virtual center} of $A$, denoted $V(A)$, is the set of
all $g \in A$ such that $[A:C_A(g)]<\infty$. This is a subgroup, because
if $g,h \in V(A)$ then the subgroup $C_A(gh)$, which contains
$C_A(g)\intersect C_A(h)$, has finite index. 

\begin{lemma}[Some characteristic subgroups]
\label{lemma:characteristic}
Let $A$ be a finitely generated, virtually abelian group. Then the
virtual center $V(A)$, its center $ZV(A)$, and its torsion subgroup
$TZV(A)$, are all characteristic for endomorphisms of $A$. Moreover,
$V(A)$ and $ZV(A)$ both have finite index in $A$, whereas $TZV(A)$ is
finite.
\end{lemma}

Lemma \ref{lemma:characteristic} is proved below.

\paragraph*{Step 6.} Consider the HNN extension $G=A_\phi$ above. Let
$V(A)$, $ZV(A)$, $TZV(A)$ be as in Lemma \ref{lemma:characteristic}, so
all these subgroups are taken into themselves by $\phi$. Since $TZV(A)$
is finite we in fact have $\phi(TZV(A)) = TZV(A)$, and so $K=TZV(A)$ is
a finite, normal subgroup of $G$.

Replacing $G$ by $G/K$, we may assume that $TZV(A)$ is trivial, and it
follows that $ZV(A)$ is torsion-free abelian, and so is isomorphic to
$\Z^n$. Since $\phi(ZV(A)) \subset ZV(A)$, the action of $\phi$ on $ZV(A)$
is given by some $n \cross n$ matrix of integers~$N$. Thus, $G/K$ has a
finite-index subgroup isomorphic to $\Gamma_{N}$, finishing the proof of
Theorem \ref{theorem:rigidity}.

\paragraph*{Proof of Lemma \ref{lemma:characteristic}.}

To see $[A:V(A)]<\infty$, if $B$ is any
finite-index abelian subgroup of $A$ then obviously $B \subset V(A)$.

Consider an endomorphism $\phi \from A \to A$. We now show that
$\phi(V(A))
\subset V(A)$. Consider $g \in V(A)$, so $[A:C_A(g)]<\infty$. It follows
that $[\phi(A):C_{\phi(A)}(\phi(g))]<\infty$, and so 
$[A:C_{\phi(A)}(\phi(g))]<\infty$. But $C_{\phi(A)}(\phi(g))
\subset C_A(\phi(g))$, and so $\phi(g) \in V(A)$.

Next we claim that $V(V(A))=V(A)$. To see why, if $g \in V(A)$ then
$[A:C_G(g)]<\infty$, and so $[V(A):C_G(g) \intersect
V(A)]<\infty$. But $C_G(g) \intersect V(A)
\subset C_{V(A)}(g)$, and so $[V(A):C_{V(A)}(g)]<\infty$,
i.e.\ $g\in V(V(A))$. 

Next we claim that $[V(A):ZV(A)]<\infty$. In fact if $V$ is any finitely
generated group which is its own virtual center, then $[V:ZV]<\infty$
(the converse is also true, trivially). To see why, let $g_1,\ldots,g_k$
be a generating set for $V$. Since $V(V)=V$, each of the groups
$C_{V}(g_1), \ldots, C_{V}(g_k)$ has finite index in $V$. It follows
that their intersection has finite index in $V$; but their intersection
is precisely $ZV$.

Now we claim that $ZV(A)$ is characteristic for endomorphisms of
$V(A)$ (and so is also characteristic for endomorphisms of $A$). In fact,
if $V$ is any finitely generated group whose center $ZV$ has finite
index, then $ZV$ is characteristic for any injective
endomorphism $\phi\from V \to V$ whose image has finite index. To see why,
we have
$Z(\phi(V)) =
\phi(ZV)$, and so
$$[\phi(V):Z(\phi(V))] = [\phi(V):\phi(ZV)] = [V:ZV]<\infty
$$
Clearly $\phi(V) \intersect ZV \subset Z(\phi(V))$, and so
$$[\phi(V):Z(\phi(V))] \le [\phi(V):\phi(V) \intersect ZV]
$$
The quotient group $V / ZV$ is finite, and the quotient homomorphism $V
\to V/ZV$, when restricted to the subgroup $\phi(V)$, has kernel $\phi(V)
\intersect ZV$. It follows that
$$[\phi(V):\phi(V) \intersect ZV] \le \abs{V/ZV} = [V:ZV] =
[\phi(V):Z(\phi(V))]
$$
All of the above inequalities are therefore equalities, and so
$$\phi(ZV) = Z(\phi(V)) = \phi(V) \intersect ZV
$$
which implies $\phi(ZV) \subset ZV$.

Finally, it is clear that for any finitely generated abelian group, the
torsion subgroup is characteristic for injective endomorphisms.

\vfill\break


\section{Questions}
\label{section:questions}

\subsection{Remarks on the polycyclic case}
\label{section:polycyclic}

Given an integer matrix $M \in \GL(n,\R)$, the group $\Gamma_M$ is
polycyclic if and only if $\abs{\det M} = 1$, and if $M \in
\GL_\cross(n,\R)$ this occurs if and only if $\Gamma_M$ is a cocompact
{\em discrete} subgroup of $G_M$. In this case it follows that $\Gamma_M$
is quasi-isometric to $G_M$, and the notion of horizontal-respecting
quasi-isometry clearly transfers to $\Gamma_M$. The techniques of this
paper do not provide a quasi-isometric classification in this case,
however they do yield the following partial result:

\begin{theorem} \quad
If $M,N \in \SL(n,\Z)$ lie on \nb{1}parameter subgroups of
$\GL(n,\R)$, then there is a horizontal  respecting quasi-isometry
$\Gamma_M
\to \Gamma_N$ if and only if there is a horizontal respecting
quasi-isometry $G_M \to G_N$, and this occurs if and only if there are
real numbers $a,b\ne 0$  such that $M^a, N^b$ have the same absolute
Jordan form.
\qed
\label{theorem:polycyclic}
\end{theorem}

This raises the question: Is every quasi-isometry $\Gamma_M \to \Gamma_N$
horizontal respecting? Equivalently, is every quasi-isometry $G_M \to G_N$
horizontal respecting? The answer is obviously no, for example when $M,N$
are identity matrices and $G_M, G_N$ are Euclidean spaces. However, we
conjecture:

\begin{conjecture}
If $M,N \in \SL(n,\Z)$ lie on \nb{1}parameter subgroups of $\GL(n,\R)$,
and if $M,N$ have no eigenvalues on the unit 
circle, then any quasi-isometry $G_M \to G_N$ is horizontal respecting.
\label{conjecture:polycyclic:horizontal}
\end{conjecture}

Moreover, Theorem \ref{theorem:polycyclic} and Conjecture
\ref{conjecture:polycyclic:horizontal} together would imply the following 
(see \cite{FarbMosher:solvable}):

\begin{conjecture}
Suppose $M\in \SL(n,\Z)$ has no eigenvalues on the unit 
circle.  If $G$ is any finitely generated group quasi-isometric to 
$\Gamma_M$, then there is a finite normal subgroup $F$ of $G$ so that 
$G/F$ is abstractly commensurable to $\Gamma_N$, for some $N \in
\SL(n,\Z)$ with no eigenvalues on the unit circle.
\label{conjecture:polycyclic:qi}
\end{conjecture}

\subsection{The quasi-isometry group of $\Gamma_M$}
\label{section:qigroup}

Given a finitely generated group $G$, the set of quasi-isometries from
$G$ to itself, modulo the identification of quasi-isometries which differ
by a bounded amount, forms a group called the \emph{quasi-isometry group}
of $G$, denoted $\QI(G)$. Given a $1 \cross 1$ matrix $M = (m)$ with $m
\ge 2$, the quasi-isometry group of the solvable Baumslag-Solitar group
$\Gamma_M \approx \BS(1,m)$ was computed in
\cite{FarbMosher:BSOne}: 
$$\QI(\BS(1,m)) \approx \Bilip(\R) \cross \Bilip(\Q_m)
$$
where $\Q_m$ is the metric space of $m$-adic rational numbers, and
$\Bilip(X)$ denotes the group of bilipschitz self maps of a metric space
$X$. 

\begin{Problem}
Compute the quasi-isometry group of $\Gamma_M$ in general.
\end{Problem}

The strongest result we have on this problem so far is
Proposition \ref{proposition:regularity}, but see the remarks after that
proposition.

In \cite{FarbMosher:BSTwo} the computation of $\QI(\BS(1,m))$ was applied
to prove quasi-isometric rigidity of $\BS(1,m)$, using techniques of
Hinkkanen \cite{Hinkkanen:quasisymmetric} and Tukia
\cite{Tukia:quasiconformal}. While quasi-isometric rigidity of $\BS(1,m)$
now has a completely different proof \cite{MosherSageevWhyte}, which we
have here generalized to $\Gamma_M$, one might still pursue:

\begin{Problem}
Give a proof of quasi-isometric rigidity of $\Gamma_M$, generalizing the
results of \cite{FarbMosher:BSTwo}.
\end{Problem}

This should lead to a deeper understanding of the geometry of
$\Gamma_M$. For example, Tukia \cite{Tukia:quasiconformal} characterizes
subgroups of the quasiconformal group of a sphere which are conjugate
into the \Mobius\ group. We have analogous results for lattices in
three-dimensional \solv-geometry, and there
should be generalizations to solvable Baumslag-Solitar groups and
to~$\Gamma_M$.




\newcommand{\etalchar}[1]{$^{#1}$}
\providecommand{\bysame}{\leavevmode\hbox to3em{\hrulefill}\thinspace}

\bigskip

\noindent
Benson Farb:\\
Department of Mathematics\\
University of Chicago\\
5734 University Ave.\\
Chicago, Il 60637\\
farb@math.uchicago.edu
\medskip

\noindent
Lee Mosher:\\
Department of Mathematics and Computer Science\\
Rutgers University, Newark\\
Newark, NJ 07102\\
mosher@andromeda.rutgers.edu

\end{document}